\newif\ifarxiv
\arxivtrue

\ifarxiv

\documentclass{article}

\else

\documentclass[aos,preprint]{imsart}

\fi

\usepackage{amsmath, amsthm, amssymb}
\usepackage{graphicx}
\usepackage{verbatim}
\usepackage{caption}
\usepackage{subcaption}
\usepackage{fancyvrb}
\usepackage{enumerate}
\usepackage{natbib}

\RequirePackage[OT1]{fontenc}
\RequirePackage[colorlinks,linkcolor=blue,citecolor=blue,urlcolor=blue]{hyperref}

\theoremstyle{plain}
\newtheorem{prop}{Proposition}

\newtheorem{coro}[prop]{Corollary}
\newtheorem{lemm}[prop]{Lemma}
\newtheorem{theo}[prop]{Theorem}

\theoremstyle{definition}

\theoremstyle{remark}

\ifarxiv
\else
\startlocaldefs
\fi

\usepackage{StefanTexCommandsV2}

\newcommand{\gn}{g^{(n)}}
\newcommand{\hgn}{\hg^{(n)}}
\newcommand{\hgnp}{\hg^{(n, \, p)}}
\newcommand{\hgns}{\hg^{(n, \, *)}}

\ifarxiv
\else
\endlocaldefs
\fi

\begin{document}

\ifarxiv

\title{The Efficiency of Density Deconvolution}
\author{Stefan Wager\footnote{I am deeply grateful to Brad Efron for many enlightening conversations as well as his continual encouragement, and to Dave Donoho and Will Fithian for several helpful comments and suggestions. This work was supported by a B. C. and E. J. Eaves Stanford Graduate Fellowship.}}

\date{Department of Statistics\\
Stanford University\\
\texttt{swager@stanford.edu}}

\maketitle

\begin{abstract}
The density deconvolution problem involves recovering a target density $g$ from a sample that has been corrupted by noise. From the perspective of Le Cam's local asymptotic normality theory, we show that \emph{non-parametric} density deconvolution with Gaussian noise behaves similarly to a low-dimensional \emph{parametric} problem that can easily be solved by maximum likelihood. This framework allows us to give a simple account of the statistical efficiency of density deconvolution and to concisely describe the effect of Gaussian noise on our ability to estimate $g$, all while relying on classical maximum likelihood theory instead of the kernel estimators typically used to study density deconvolution.

\end{abstract}

\fi

\section{Introduction}

Suppose that we observe $n$ samples $X_i$ drawn from a hierarchical model
\begin{align}
\label{eq:setup}
\mu_i \simiid g\p{\cdot}, \ \
X_i = \mu_i + \varepsilon_i, \ \ \text{with}  \ \
\varepsilon_i \simiid \nn\p{0, \, 1},
\end{align}
and our goal is to estimate the unknown density $g(\cdot)$. This problem, sometimes called the density deconvolution problem, is remarkably hard in terms of asymptotic statistical criteria. For example, as shown by \citet{carroll1988optimal} and \citet{fan1991optimal}, if we assume a non-parametric setup where $g$ is only known to have a Lipschitz-continuous $k$-th derivative, then the minimax error rate for estimating $g$ under the integrated squared error loss decays as $\log(n)^{-(k + 1)}$.

The density deconvolution problem is traditionally studied using kernel methods. The motivation for kernel estimators is that they often achieve adaptive minimaxity over natural regularity classes for $g$; for example, they in fact achieve the optimal rates of \citet{carroll1988optimal} and \citet{fan1991optimal} described above. The properties of kernel density deconvolution have been analyzed by several authors, including \citet{butucea2009adaptive}, \citet{carroll2004low}, \citet{comte2011data}, \citet{efromovich1997density}, \citet{fan2002wavelet}, \citet{hall2008estimation}, \citet{hall2007ridge}, \citet{stefanski1990deconvoluting}, \citet{wand1998finite}, and \citet{zhang1990fourier}.

Despite the prolific literature devoted to them, however, kernel methods do not necessarily yield a fully satisfying theory of the statistics of density deconvolution. In particular, \citet{efron2014bayes,efron2014two} proposed a simple maximum likelihood approach to density deconvolution that performs qualitatively better than kernel methods on several realistic scientific tasks. Thus, it appears that while kernel methods are nearly optimal in terms of the standard asymptotic criteria used in the literature, they are not always optimal in practical applications. This suggests a need for a new optimality theory for density deconvolution.

The goal of this paper is to move us towards such a theory. We begin with a close analysis of Efron's method, which involves estimating the unknown density $g$ by maximum likelihood with a $p$-parameter model
\begin{equation}
\label{eq:expfam}
g_\eta\p{\mu} = g_0\p{\mu} \, \exp\sqb{\eta \cdot T\p{\mu} - \psi\p{\eta}},
\end{equation}
where $T\p{\mu}$ is some carefully chosen $p$-parameter statistic and $\psi\p{\cdot}$ is the log-partition function. Given that Efron's method involves \emph{parametric} maximum likelihood estimation, it may appear surprising that this method would work well in a \emph{non-parametric} setup where we only know that the target density $g$ belongs to some regularity class. However,  we show that maximum likelihood estimation in the model \eqref{eq:expfam} with an appropriate choice of $T$ has adaptive minimaxity properties that are reminiscent of those enjoyed by kernel estimators. Moreover, because Efron's method relies on maximum likelihood instead of the ad-hoc kernel inversion procedure used by classical methods, we may also hope for the method to be well-behaved in a wide variety of practical applications.

We build our analysis around a \emph{relative efficiency} criterion, which measures the information loss for estimation in the model \eqref{eq:expfam} we incur from the noise $\varepsilon$ in \eqref{eq:setup}. Our main result is that, given a specified carrier $g_0$, there exists a low-dimensional statistic $T$ that captures essentially all the information contained in the $X$-sample for estimating local perturbations to $g_0$. Any model of the form \eqref{eq:expfam} that tries to use a higher-dimensional parametrization will have some bad contrasts that could have been accurately estimated from clean observations $\mu_i$, but become effectively impossible to estimate from the noised observations $X_i$. If $g_0$ is Gaussian, the number of samples needed to accurately estimate a $p$-parameter family of the form \eqref{eq:expfam} scales exponentially in $p$, and this bound holds uniformly over any possible choice of the statistic $T$. 

At face value, our results can be interpreted as a companion to the classical hardness results of \citet{carroll1988optimal} and \citet{fan1991optimal}: our analysis implies that it is impossible to accurately estimate any parametric family of the form \eqref{eq:expfam} in terms of relative efficiency if $p$ is even moderately large. Thus, no matter how clever we may be, we cannot use Efron's density deconvolution model to efficiently learn a rich model for $g$.

However, from a decision-theoretic point of view, our result can also be interpreted in a more optimistic light. Because high-dimensional models of the form \eqref{eq:expfam} are effectively impossible to estimate, we lose almost nothing by just using a low dimensional model. Thus, the non-parametric density deconvolution problem effectively reduces to parametric inference for Efron's model \eqref{eq:expfam} with an appropriate low-dimensional parametrization $T$.

Our results take on a particularly simple form when the carrier $g_0$ is Gaussian. In this case, the ``optimal'' statistics are polynomials in $\mu$, yielding the class of estimators
\begin{equation}
\label{eq:poly}
\hg_\eta^p\p{\mu} = \exp\sqb{\sum_{j = 1}^p \eta_j \, \mu^j - \psi_p\p{\eta}},
\end{equation}
where $p$ is a tuning parameter.
As we will show, given an appropriate (usually small) choice of $p$, parametric inference for this deceptively simple model is nearly equivalent to optimal non-parametric inference for $g$.\footnote{In all our experiments, we use $p = 4$. In practice, $p$ could also be selected by cross-validation.} More specifically, this class of estimators is asymptotically nearly adaptively minimax under Kullback-Leibler (KL) loss for estimating local perturbations of $g_0$
\begin{equation}
\label{eq:local}
g\p{\mu} = g_0\p{\mu} \exp\sqb{\tau\p{\mu} - \psi_\tau}, \ \ g_0\p{\mu} = \frac{1}{\sigma} \, \varphi\p{\frac{\mu}{\sigma}}, 
\end{equation}
where $\tau(\cdot)$ is a tilting function satisfying regularity conditions detailed in Section \ref{sec:minimax}. Moreover, the efficiency shortfall of the best estimator of the form \eqref{eq:poly} relative to the minimax estimator can be bounded by a small explicit constant that is, for typical parameter values, on the order of 2.

In summary, this paper introduces a relative efficiency criterion and a local perturbation model that, together, enable us to shed new light on the classical problem of density deconvolution. We show that the exponential family method \eqref{eq:expfam} cannot get around standard non-parametric impossibility results and, in particular, does not allow us to estimate rich models for $g(\cdot)$. However, the method does an excellent job at extracting all the available information from the $X$-sample, all while allowing for straight-forward estimation and inference. Thus, despite its simple form, Efron's parametric density deconvolution method yields simple estimators whose excellent practical performance is firmly grounded in asymptotic minimaxity theory.

\subsection{Theoretical Setup: Local Deconvolution}
\label{sec:setup}

Throughout this paper, we assume that
we observe $n$ samples $X_i$ drawn from a hierarchical model
\begin{align}
\label{eq:setup2}
&\mu_i \simiid g\p{\cdot}, \ 
X_i = \mu_i + \varepsilon_i, \ 
\varepsilon_i \simiid \nn\p{0, \, 1}; \text{ thus,} \\
\label{eq:fdef}
&X_i \simiid f\p{\cdot}, \ f\p{x} = \p{\varphi * g}\p{x},
\end{align}
where $f$ denotes the marginal density of the observations $X$ generated according to the model \eqref{eq:setup} and $\varphi$ is the standard Gaussian density. The setup induced by \eqref{eq:setup2} and \eqref{eq:fdef} is classical; however, the analysis techniques we use to understand this model are rather different from the ones used by, e.g., \citet{carroll1988optimal}, \citet{efromovich1997density}, or \citet{fan1991optimal}.

Unlike other authors who model $g(\cdot)$ as a fixed member of a regularity class $\rr$ that can be estimated with increasing accuracy as the sample size $n$ grows to infinity, we study a sequence of ever-shrinking perturbations of a known carrier density $g_0$:
\begin{equation}
\label{eq:perturb_setup}
\gn \p{\mu} = g_0\p{\mu} \exp\sqb{\frac{1}{\sqrt{n}} \tau\p{\mu} - \psi_n}, 
\end{equation}
where $\tau\p{\mu}$ belongs to an appropriate regularity ellipsoid. Notice that, as $n$ gets large, estimating $\tau$ does not necessarily get easier because the deviance between $g_0$ and $\gn$ decays as $1/\sqrt{n}$. As we will show formally in Section \ref{sec:minimax}, the problem of estimating $\gn$ in \eqref{eq:perturb_setup} under a re-scaled deviance loss
\begin{equation}
\label{eq:deviance}
L_n\p{\hgn} = n\, D_{KL}\p{\gn, \, \hgn} = n \int \gn\p{\mu} \log\p{\frac{\gn\p{\mu}}{\hgn\p{\mu}}} \, d\mu.
\end{equation}
converges to a \emph{locally asymptotically normal} experiment in the sense of \citet{lecam1960locally}. The formalism provided by the sequence of problems \eqref{eq:perturb_setup} thus enables us to study density deconvolution from the perspective of classical maximum likelihood asymptotics.

Although the estimation problem in the local perturbation model \eqref{eq:perturb_setup} may look very different from the classical estimation problem with $g \in \rr$ for some fixed regularity set $\rr$, it appears that studying the former can help us learn about the latter. First of all, we find a universality phenomenon, where maximum likelihood estimation in the model \eqref{eq:poly} is simultaneously nearly minimax for any Gaussian carrier $g_0 = \varphi_\sigma$, and so we can carry out practical data analysis with a simple default model for $g$. Second, our analysis recovers familiar qualitative aspects of the theory of \citet{carroll1988optimal} and \citet{fan1991optimal}, such as exponential blow-ups in the number of samples required to estimate more complex models; see, e.g., Theorem \ref{theo:gauss_hard}. Thus, insofar as maximum-likelihood estimation in our local perturbation model is amenable to exact asymptotic analysis, it appears that the model \eqref{eq:perturb_setup} is a useful theoretical tool for gaining new insights about the density deconvolution problem.

We begin our analysis by developing a theory of relative efficiency for density deconvolution in Section \ref{sec:efficiency}, with the goal of establishing lower bounds for the error rate of any $p$-parameter model of the form \eqref{eq:expfam}. In Section \ref{sec:minimax}, we complement this worst-case picture by establishing near-minimax properties of the estimator \eqref{eq:poly} for estimating local perturbations of Gaussian densities. Our proof technique is built on the local asymptotic normality theory of \citet{lecam1960locally} combined with Pinsker's theorem for the Gaussian sequence model. Finally, in Section \ref{sec:examples}, we discuss practical differences between kernel density deconvolution and Efron's method, and the respective optimality theories that justify each method. All proofs are provided in the appendix.

Our minimaxity proof, which first reduces a continuous problem to a Gaussian sequence model and then applies Pinsker's theorem,
fits into a rich literature on solving non-parametric problems via Gaussian estimation; see, e.g., \citet{brown2004equivalence}, \citet{brown1996asymptotic}, \citet{efromovich1996asymptotic}, \citet{golubev2010asymptotic}, \citet{johnstone2011gaussian}, and \citet{nussbaum1996asymptotic}.
The problem of density estimation using a log-spline model for $g$ of the form \eqref{eq:expfam} has also been considered by \citet{koo1999logspline} and \citet{koo1998log}; however, their approach more closely follows that traditional ideas of \citet{fan1991optimal} and others.
Finally, we note that the non-parametric maximum likelihood problem of estimating $f = \varphi * g$ has been studied by, among others, \citet{jiang2009general}, \citet{koenker2014convex}, and \citet{zhang2009generalized}. Efron's method also induces a natural estimator $f_{\heta} = \varphi * g_{\heta}$, but we do not study its properties here.

\paragraph{Remark: Cyclic Convolution}

To avoid technical difficulties relating to integrability over unbounded domains, we follow convention and study a \emph{cyclic} convolution model \citep[e.g.,][]{efromovich1997density}. More specifically, we assume that our observations are within a bounded interval $\Omega_M = [-M, \, M]$, and that we have a cyclic convolution operator
\begin{equation}
\label{eq:KM}
K_M : \Omega_M \times \Omega_M \rightarrow \RR_+, \ \ K_M\p{\mu, \, x} = \sum_{j = -\infty}^{\infty} \varphi\p{x - \mu + 2 jM}.
\end{equation}
We also frequently use the cyclically wrapped-around Gaussian density with variance $\sigma^2$, i.e., $\smash{\varphi^M_\sigma(x) = \sum_{j = -\infty}^{\infty} \sigma^{-1} \, \varphi\p{\p{x + 2 jM} / \sigma}}$ with $x \in \Omega_M$, as the carrier $g_0$.
In our results, $M$ should be thought of as large enough that $K_M\p{\mu, \, x} \approx \varphi\p{x - \mu}$ over the range of the data; formally, we will seek to state results in a limit with $M \rightarrow \infty$.

\section{Relative Efficiency and Most Stable Families}
\label{sec:efficiency}

The observations $X_i$ generated by the model \eqref{eq:setup} are noisy measurements of clean draws $\mu_i$ from $g$. What is the information loss due to this extra noise? Or, in other words, if we can estimate $\eta$ to a given accuracy using $n_\mu$ samples $\mu$ drawn directly from $g$, how many samples $n_X$ would we need to estimate $\eta$ to the same accuracy using only $X$ samples?
Although this question may not at first appear directly related to our topic of interest, answering it will prove to be helpful in guiding us towards good choices for the statistic $T$ in Efron's model \eqref{eq:expfam} and in proving lower bounds for the error rate of any parametric model.

In Section \ref{sec:relative} below, we define a {\it relative efficiency coefficient} that will let us make the above question precise. Using this formalism, we then derive the ``best'' families of the form \eqref{eq:expfam} in terms of relative efficiency in Section \ref{sec:stable}. Finally, in Sections \ref{sec:hardness} and \ref{sec:nongauss} we interpret these relative efficiency results for both Gaussian and non-Gaussian carriers $g_0$ respectively.

\subsection{The Relative Efficiency Coefficient}
\label{sec:relative}

If we had access to samples $\mu_i$ drawn directly from a density $g_\eta$ of the form \eqref{eq:expfam}, then under standard regularity conditions maximum likelihood estimation would yield an asymptotically normal estimator $\heta_\mu$ with
\begin{equation}
\sqrt{n} \p{\heta_\mu - \eta} \rightarrow \nn\p{0, \, \ii_\mu^{-1}\p{\eta}}, \ \ \ii_\mu\p{\eta} = -\EE[\eta]{\frac{\partial^2}{\partial \eta^2} \log g_\eta\p{\mu}},
\end{equation}
where $\ii_\mu\p{\eta}$ is the Fisher information for estimating $\eta$ carried by the $\mu$-samples.
In our setup, however, we only have access to samples $X$ drawn from \eqref{eq:fdef};
the Fisher information for the maximum-likelihood estimator $\heta_X$ is then reduced to
\begin{equation}
\label{eq:fi}
\ii_X\p{\eta} = -\EE[\eta]{\frac{\partial^2}{\partial \eta^2} \log f_\eta\p{X}}.
\end{equation}
Given this setup, we define the efficiency of $\heta_X$ relative to $\heta_\mu$ as
\begin{equation}
\label{eq:rel_mult}
\rho_\eta\p{T} = \inf_{a \in \RR^p, \, a \neq 0} \cb{\frac{a^\top \ii_X \p{\eta} a}{a^\top \ii_\mu\p{\eta} a}},
\end{equation}
and argue that the relative efficiency coefficient $\rho$ provides a natural measure of information loss due to the noise $\varepsilon$ in \eqref{eq:setup}.

If we wanted to measure relative efficiency in a univariate family \eqref{eq:expfam} with tilting function $t: \RR \rightarrow \RR$, then the natural measure of relative efficiency is
\begin{equation}
\label{eq:relative}
\rho_{\eta}\p{t} = \limn \frac{\Var{\heta_\mu}}{\Var{\heta_X}} = \frac{\ii_X\p{\eta}}{\ii_\mu\p{\eta}},
\end{equation}
which measures exactly the increase in sample size required for accurate estimation using the $X$-sample instead of the $\mu$-sample. As we can easily verify, the multivariate relative information coefficient \eqref{eq:rel_mult} is simply the worst-case relative efficiency for any univariate subfamily of \eqref{eq:expfam}:
\begin{equation}
\label{eq:rel_mult2}
\rho_{\eta}\p{T} = \inf\cb{\rho_\eta\p{t} : t\p{\mu} = a \cdot T\p{\mu}, \, a \in \RR^p}. 
\end{equation}
This connection provides a first motivation for the definition \eqref{eq:rel_mult}.
The coefficient $\rho_\eta$ can also be viewed as a natural extension to the classical $E$-optimal criterion for experimental design \citep{ehrenfeld1955efficiency}:  if $T$ is scaled such that $\ii_\mu = I_{p \times p}$, then $\rho_\eta\p{T}$ is the minimum eigenvalue of $\ii_X$---which is exactly the criterion that the $E$-optimal designs maximize.

The following result derives a simple functional form for the relative efficiency coefficient. We note a striking similarity between the formula \eqref{eq:variance} and the results of \citet{louis1982finding} for maximum likelihood estimation with missing data. This similarity is not an accident: we could also interpret the density deconvolution problem as a one where we would would have wanted to observe $\mu = X - \varepsilon$, but $\varepsilon$ is missing.

We also establish a key theoretical property of $\rho$, namely that it is transformation invariant: if we apply any invertible linear transformation $Q$ to $T$, the value of $\rho$ remains unchanged. This transformation invariance provides further evidence that $\rho$ is a ``natural'' measure for understanding the difficulty of density deconvolution. We note in particular that the Fisher information $\ii_X$ is not transformation invariant. 

\begin{lemm}
\label{lemm:relative}
The relative efficiency coefficient \eqref{eq:rel_mult} is transformation invariant, in the sense that $\rho_{0}\p{T} = \rho_{0}\p{QT}$ for any invertible linear transformation $Q$. For a univariate statistic, the relative efficiency coefficient can be written as
\begin{equation}
\label{eq:variance}
\rho_{\eta}\p{t} = \frac{\Var[\eta]{\EE{t\p{\mu} \cond X}}}{\Var[\eta]{t\p{\mu}}}.
\end{equation}
\end{lemm}

\subsection{Deriving the Most Stable $p$-dimensional Family}
\label{sec:stable}

Given our notion of relative efficiency defined above, it is natural to ask whether there exist optimal $p$-dimensional statistics $T$ in terms of this criterion. Perhaps surprisingly, we will show that not only do such optimal statistics exist, but they are in general quite easy to compute and can give us guidance for practical data analysis. Formally, we define the optimal statistics as solutions to the optimization problem
\begin{align}
\label{eq:gamma}
\Gamma_p\p{g_0} = \argmax_{T : \Omega \rightarrow \RR^P} \cb{\rho_0\p{T}},
\end{align}
where $g_\eta\p{\mu} = g_0\p{\mu} \exp\sqb{\eta \cdot T\p{\mu} - \psi\p{\mu}}$ is defined in terms of some known carrier.  Notice that our definition of $\Gamma$ in terms of relative efficiency at $\eta = 0$ is without loss of generality, since we could always just use $g_\eta(\mu)$ as our carrier if this condition did not hold.

In order to solve the problem \eqref{eq:gamma}, we begin with a technical result for computing the multivariate relative efficiency coefficient. In this section, we will assume that the $X$ and $\mu$ are distributed over a compact interval $\Omega$, and that $X$ is noised by a generic convolution operator $K\p{\mu, \, x}$.

\begin{lemm}
\label{lemm:mult}
Let $T$ be a $p$-dimensional statistic; let $g\p{\cdot} = g_\eta\p{\cdot}$ be defined as in \eqref{eq:expfam} and let $f$ be the marginal density of the observations $X$. Then
\begin{align}
\label{eq:mult}
&\rho_\eta\p{T} = \rho_\eta\p{a^* \cdot T}, \with\\
\label{eq:mult2}
&a^* = \underset{\Norm{a}_2 = 1}{\argmin} \cb{\frac{a^\top \int_\Omega T\p{\mu} T^\top\!\!\p{\mu} K^2\p{x, \, \mu} g^2\p{\mu}f^{-1}\p{x}  dx \, d\mu \ a}{a^\top \int_\Omega T\p{\mu} T^\top\!\!\p{\mu} g\p{\mu} d\mu \ a}}.
\end{align}
\end{lemm}

We note that, in practice, the optimization problem described in \eqref{eq:mult2} below can be efficiently solved by taking $\smash{a^* = Q_T^{-1} b^* / ||Q_T^{-1} b^*||_2}$, where $b^*$ is the eigenvector corresponding to the smallest eigenvalue of $Q_T^\top \, M_T \, Q_T$, and the matrices $M_T$ and $Q_T$ are defined by
\begin{align}
\label{eq:prag1}
&M_T =  \int_\Omega T\p{\mu} T^\top\!\!\p{\mu} K^2\p{x, \, \mu} g^2\p{\mu}f^{-1}\p{x}  dx \, d\mu, \\
\label{eq:prag2}
&Q_T^\top \ \int_\Omega T\p{\mu} T^\top\!\!\p{\mu} g\p{\mu} d\mu \ Q_T = I_{p \times p}.
\end{align}
With this result in hand, we are now ready to find the most favorable statistic $\Gamma_p(g_0)$. Our construction hinges around eigenfunctions of the linear operator
\begin{equation}
\label{eq:lin}
P_g\p{\mu_1, \, \mu_2} = \int_\Omega \sqrt{g\p{\mu_1}} K\p{\mu_1, \, x} f^{-1}\p{x} K\p{x, \, \mu_2} \sqrt{g\p{\mu_2}} \ dx,
\end{equation}
where $f$ is as usual defined as $f = \varphi * g$. We can verify that the top eigenfunction of $P_g$ is given by $\sqrt{g\p{\cdot}}$ and has eigenvalue 1; the $p$ subsequent eigenvectors then generate the most favorable $p$-dimensional family for density deconvolution.

\begin{theo}
\label{theo:most_favorable}
Suppose that both the carrier density $g_0 : \Omega \rightarrow \RR_+$ and the kernel $K: \Omega^2 \rightarrow \RR_+$ are continuous and bounded away from $0$, and that $P_{g_0}$ as defined in \eqref{eq:lin} is a compact operator over the space $L_2\p{\Omega}$ of square-integrable functions over $\Omega$. Then, $P_{g_0}$ admits a spectral decomposition $\zeta_1, \, \zeta_2, \, ...$, where the first eigenfunction $\zeta_1\p{\mu} = \sqrt{g_0\p{\mu}}$ has an eigenvalue 1. Moreover, the $p$-dimensional exponential family of the form \eqref{eq:expfam} with statistics
\begin{equation}
T_j\p{\mu} = \frac{1}{\sqrt{g_0\p{\mu}}} \, \zeta_{j+1}\p{\mu}
\end{equation}
is a most favorable $p$-dimensional family in the sense of \eqref{eq:gamma}, and the relative efficiency coefficient corresponds to the $p+1$-st eigenvalue of $P_{g_0}$. Finally, if the spectrum of $P_{g_0}$ does not have repeated eigenvalues, the most favorable family $\Gamma_p\p{g_0}$ is unique up to scaling and rotation.
\end{theo}

In other words, Theorem \ref{theo:most_favorable} establishes a direct link between the difficulty of learning rich perturbation models around $g_0$, and the decay rate of the spectrum of the linear operator $P_{g_0}$: the best-case efficiency for learning a $p$-parameter model depends on the $p$-th non-trivial eigenvalue of $P_{g_0}$.

As a corollary to this result, we also see that if the relative information for learning the most favorable $p$-parameter family is small, then the model \eqref{eq:expfam} with statistics $\Gamma_p(g_0)$ provides a good $X$-space approximation to any local perturbation of $g_0$. Recall that $\smash{D_{KL}(f, \, f')}$ effectively measures the power of the $X$-sample likelihood-ratio test for distinguishing $f$ from $f'$; thus, the result below implies that if $\smash{\lambda_{p + 2} \p{P_{g_0}}}$ is small, it is statistically impossible to detect any deviations from the most favorable $p$-parameter family using only $X$-samples.

\begin{coro}
\label{coro:approx}
Suppose we have a data-generating function of the form
\begin{equation}
\mu_1, \, ..., \, \mu_n \sim g^{(n)}_\tau\p{\mu} = g_0\p{\mu} \exp\sqb{\frac{\tau\p{\mu}}{\sqrt{n}} - \psi\p{\frac{\tau}{\sqrt{n}}}}
\end{equation}
for some square-integrable function $\tau$ satisfying $\smash{\int_\RR \tau^2\p{\mu} g_0\p{\mu} \ d\mu \leq C^2}$, and write $\smash{f^{(n)}_\tau = \phi * g^{(n)}_\tau}$.
Then, under the conditions of Theorem \ref{theo:most_favorable}, there exists another tilting function $\smash{\tau^{(p, \, *)}}$ in the span of the most-favorable $p$-dimensional family $\smash{\Gamma_p\p{g_0}}$ that can closely approximate $\smash{f^{(n)}_\tau}$:
\begin{equation}
\label{eq:floss}
\tau^{(p, \, *)} = \sum_{j = 1}^p \gamma_j^* \, T_j\p{\mu}, \ 
\limn n \, D_{KL} \p{f^{(n)}_\tau, \, f^{(n)}_{\tau^{(p, \, *)}}} \leq \frac{1}{2} \, C^2 \, \lambda_{p + 2} \p{P_{g_0}},
\end{equation}
where $\smash{\lambda_{p + 2} \p{P_{g_0}}}$ denotes the $p+2$-nd eigenvalue of the linear operator defined in \eqref{eq:lin}.
\end{coro}

\subsection{The Hardness of Local Deconvolution near Gaussian Carriers}
\label{sec:hardness}

The relative efficiency bound given in Theorem \ref{theo:most_favorable} is quite general, but understanding its implications for practical data analysis may not be trivial at first glance. If we are willing to take the carrier density $g_0$ to be Gaussian with variance $\sigma^2$, however, the quantities defined in the statement of Theorem \ref{theo:most_favorable} take on simple and interpretable forms. In the result below, we assume that we are doing cyclic deconvolution over the domain $\Omega_M = [-M, \, M]$ for some large $M$; this lets us state a clean result while avoiding integrability concerns over unbounded domains.

\begin{theo}
\label{theo:gauss_hard}
Suppose that the carrier $g_0\p{\mu} = \varphi^M_\sigma\p{\mu}$ is the Gaussian with variance $\sigma^2$ wrapped over the interval $\Omega_M = [-M, \, M]$, and that we are in the cyclic convolution setup over $\Omega_M$ described in Section \ref{sec:setup}. Then, for any $p$-parameter statistic $T$,
\begin{equation}
\label{eq:gauss_hard}
\rho_0\p{T} \leq \frac{1}{\p{1 + \sigma^{-2}}^p} + o_M(1),
\end{equation}
where the residual term $o_M(1)$ becomes negligible as $M$ gets large.
Moreover, this bound is satisfied by the model \eqref{eq:expfam} whose statistics are the normalized Hermite polynomials
\begin{equation}
\label{eq:hermite}
\Gamma_j\p{\mu} = \frac{1}{\sigma} H_j\p{\frac{\mu}{\sigma}}, \ \ H_j\p{\mu} = \frac{1}{\sqrt{j!}} \exp\p{\frac{\mu^2}{2}} \p{\frac{\partial}{\partial \mu}}^j \exp\p{-\frac{\mu^2}{2}},
\end{equation}
for which $\ii_\mu = I_{p \times p}$, and $\ii_X$ is diagonal with entries $\p{1 + \sigma^{-2}}^{-1}$, ..., $\p{1 + \sigma^{-2}}^{-p}$.
\end{theo}

In other words, Theorem \ref{theo:gauss_hard} implies that, if we are trying to solve the problem \eqref{eq:expfam} with a $p$-dimensional statistic $T$ scaled to have $\ii_\mu(0) = I_{p \times p}$, then we are necessarily faced with a difficult 1-dimensional sub-family with information bounded by $(1 + \sigma^{-2})^{-p}$. In particular, it is impossible to accurately distinguish local alternatives of $0$ for $\eta$ with less than $n \sim (1 + \sigma^{-2})^{p}$ observations, i.e., regardless of our choice of $T$, the number of samples required for accurate estimation scales exponentially in $p$.

Qualitatively, this exponential scaling in $p$ provides a direct analogue to the logarithmic convergence rates for kernel density estimation derived by \citet{carroll1988optimal} and \citet{fan1991optimal}:
If we work in an asymptotic regime where we incur a bias on the order of $p^{-\zeta}$ from learning a model with only $p$ degrees of freedom, then our result suggests that---heuristically---the error rate cannot decay faster than $\smash{\log_{1 + \sigma^{-2}}(n)^{-\zeta}}$. The reason this analogy is only heuristic is that our theory uses a local perturbation model, whereas that of \citet{carroll1988optimal} and \citet{fan1991optimal} is global. The connection between the two theories is however encouraging, in that it suggests that our local perturbation model enables us to get a nuanced grasp of the statistics on density deconvolution without changing the fundamental nature of the problem.

An important consequence of Theorem \ref{theo:gauss_hard} is that, for all $p \geq 2$ and any $\sigma > 0$, the family of distributions attaining the bound \eqref{eq:gauss_hard}, namely
\begin{equation}
\label{eq:hermite_fam}
g_\eta\p{\mu} = \frac{1}{\sigma} \, \varphi\p{\frac{\mu}{\sigma}} \exp \sqb{\sum_{j = 1}^p \eta_j H_j\p{\frac{\mu}{\sigma}} - \psi\p{\eta}},
\end{equation}
is equivalent to the family \eqref{eq:poly} after re-parametrizing $\eta$. In other words, we find that the polynomial log-density model \eqref{eq:poly} is universally the most-favorable family for density deconvolution near Gaussian carriers. Moreover, as a direct consequence of Corollary \ref{coro:approx}, we see that this family can be used to closely approximate any local perturbation to a Gaussian.

In Section \ref{sec:minimax}, we will provide further justification for the estimator \eqref{eq:hermite_fam} by establishing minimax-optimality properties. Before doing so, however, we pause to study the implications of Theorem \ref{theo:most_favorable} in the case when $g_0$ is not Gaussian.

\subsection{Local Density Deconvolution with Non-Gaussian Carriers}
\label{sec:nongauss}

We can also use Theorem \ref{theo:most_favorable} to compute most favorable families for estimating perturbations to non-Gaussian carriers $g_0$. In this case, we can no longer derive closed-form expressions for them; however, we can still proceed numerically.

\begin{figure}[p]
\centering
\begin{tabular}{cc}
\includegraphics[width=0.5\textwidth]{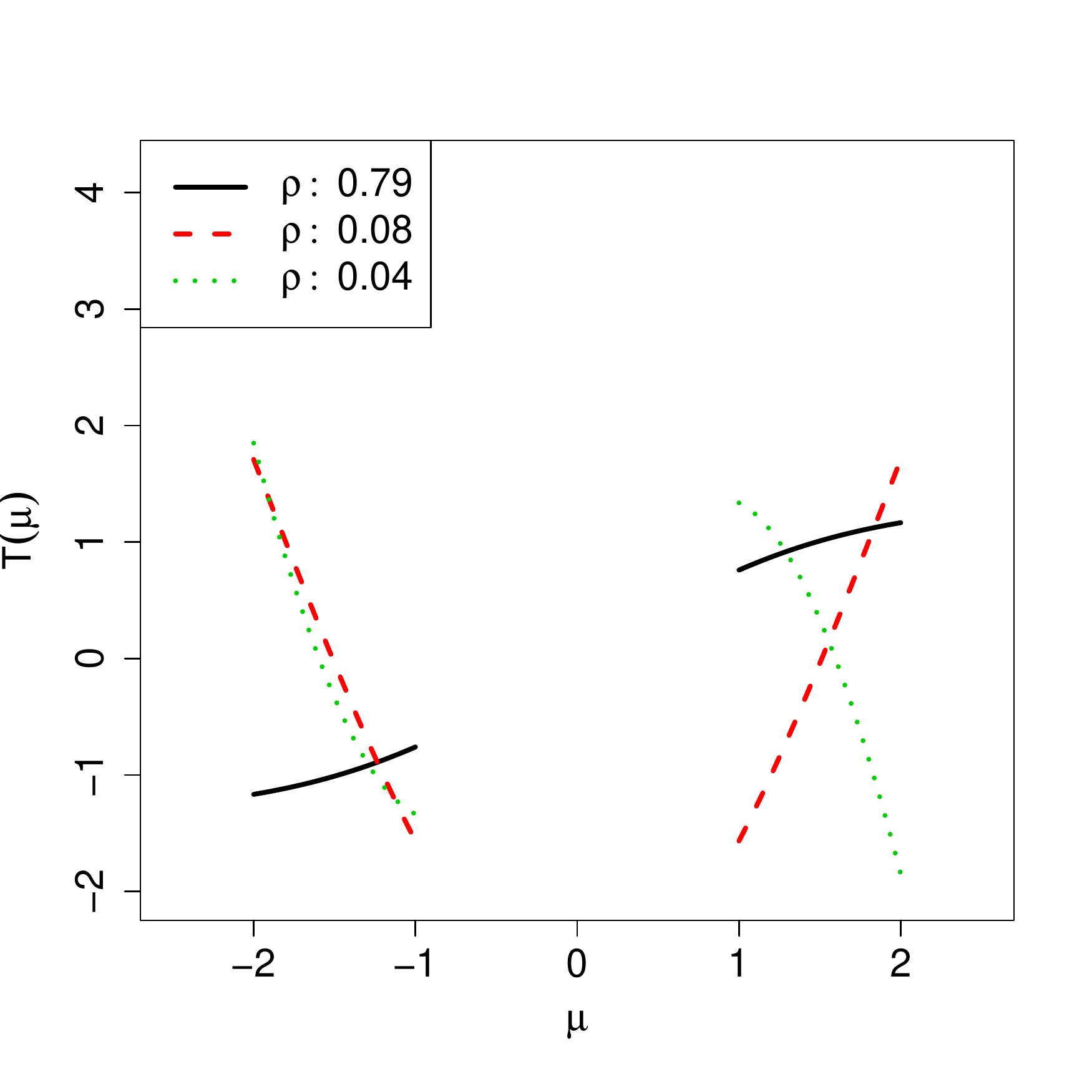} &
\includegraphics[width=0.5\textwidth]{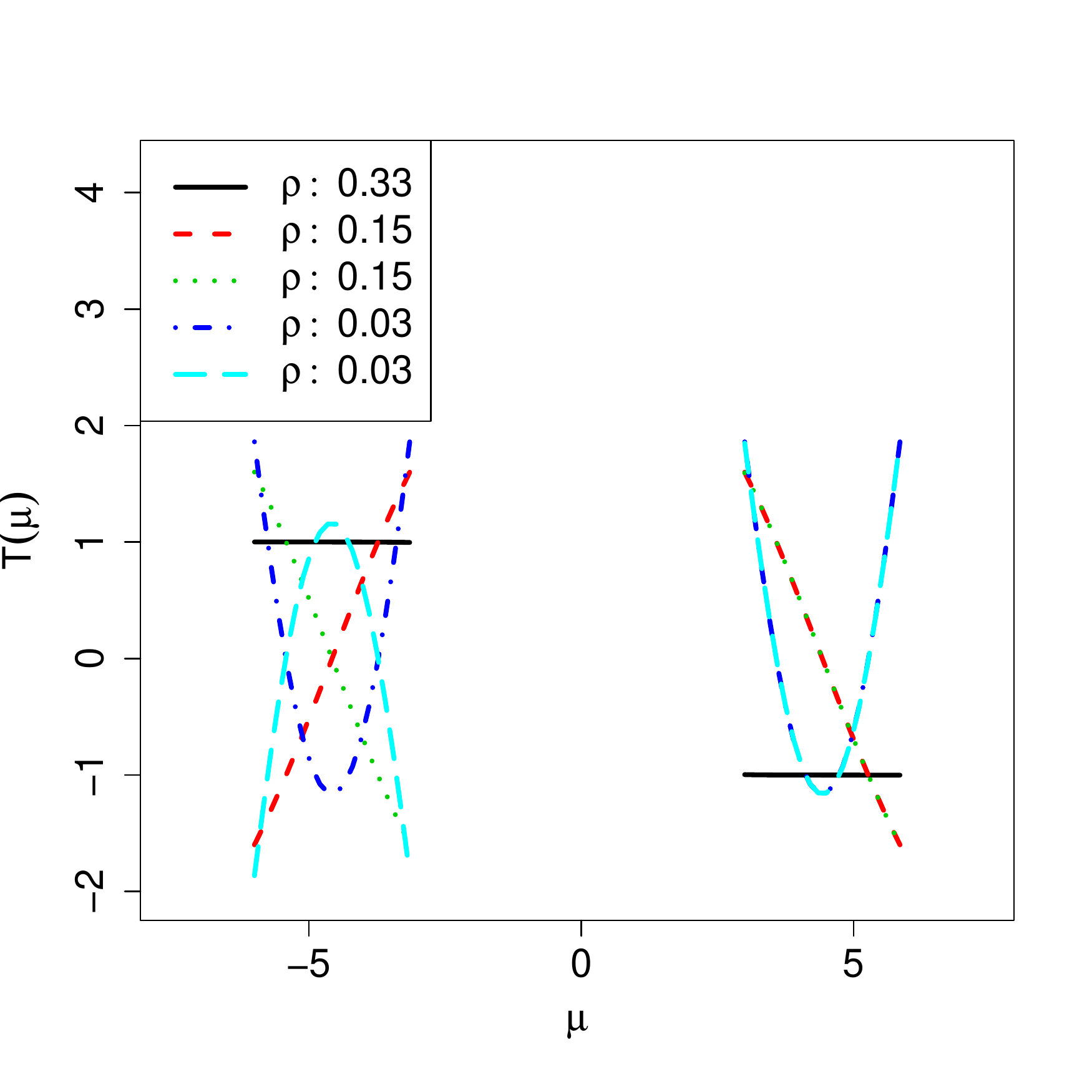} \\
Two narrow towers & Two thick towers \\
\includegraphics[width=0.5\textwidth]{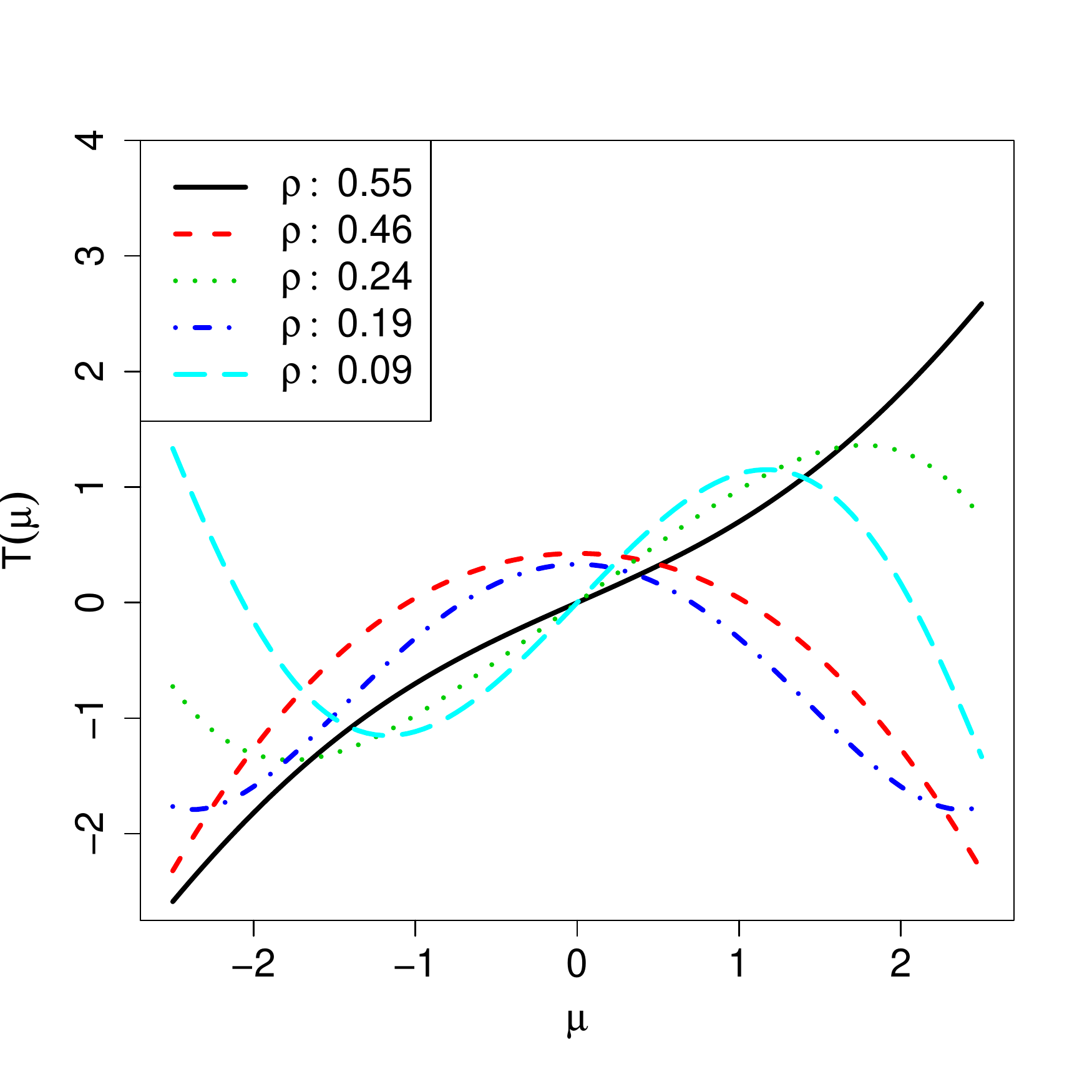} &
\includegraphics[width=0.5\textwidth]{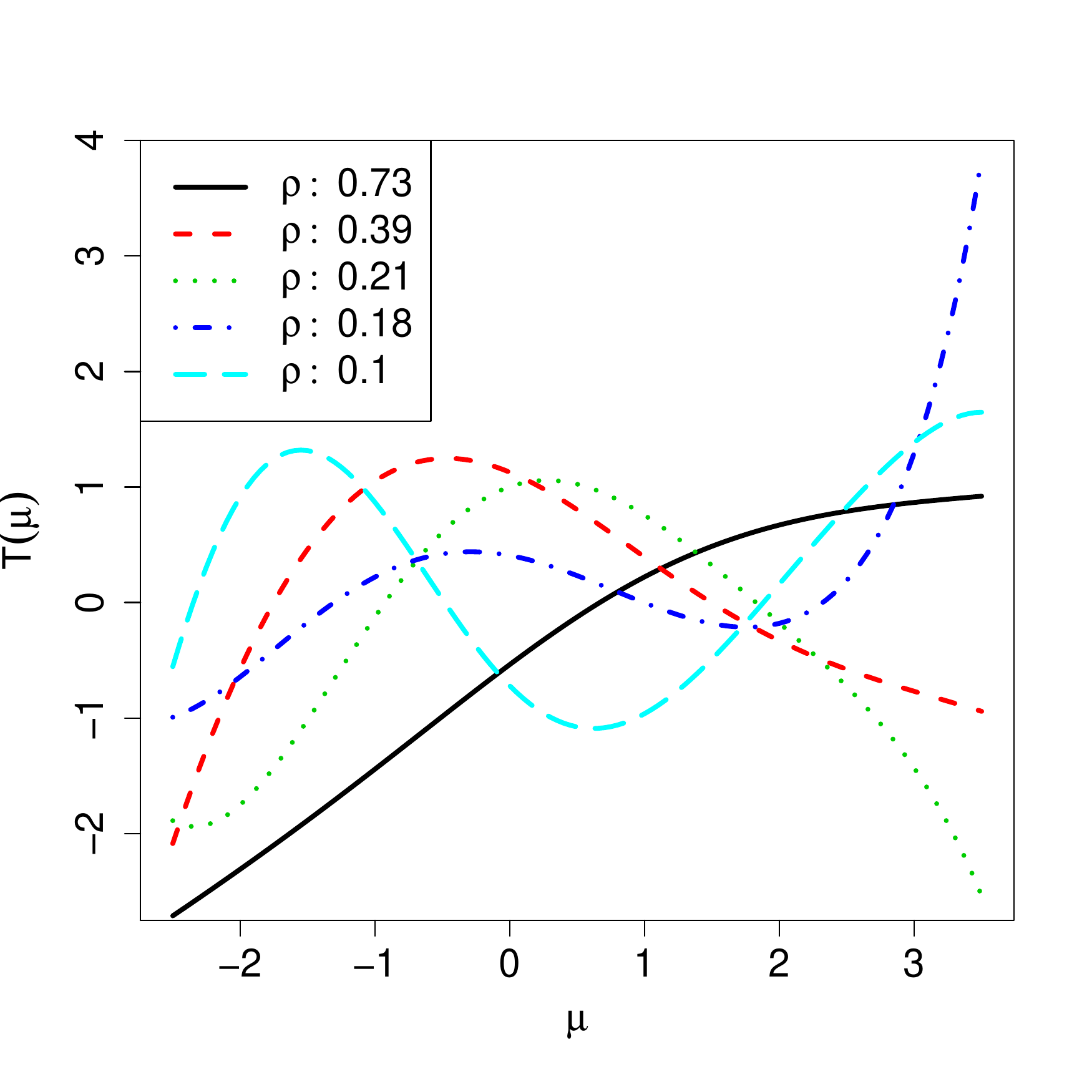} \\
Gaussian with Spike at 0 & Gaussian with Spike at 2 \\
\end{tabular}
\caption{Most favorable statistics for various carriers non-Gaussian carriers $g_0$. In each case, the noise was standard Gaussian $\varepsilon \sim \nn\p{0, \, 1}$; the different carriers $g_0$ are described in Section \ref{sec:nongauss}. The relative information coefficient $\rho$ is computed univariately for each candidate statistic.}
\label{fig:nongauss}
\end{figure}

We begin by examining a ``two towers'' model, also considered by \citet{efron2014two}. We vary the scale of the carrier
$$ g_0^{(TT, \, 1)}\p{\mu} = \frac{1}{2} \, 1\p{\cb{1 \leq \abs{\mu} \leq 2}} \  \eqand \  g_0^{(TT, \, 2)}\p{\mu} = \frac{1}{6} \, 1\p{\cb{3 \leq \abs{\mu} \leq 6}}. $$
The most favorable statistics for each choice are shown in the top row of Figure \ref{fig:nongauss}. Interestingly, for the thin tower model $g_0^{(TT, \, 1)}$, the towers appear to be too close to each other for them to be adequately distinguishable after adding standard Gaussian noise to the clean observations $\mu_i$; thus, the most favorable statistic is effectively the mean $T\p{\mu} = \mu$. Conversely, with $g_0^{(TT, \, 2)}$, the most favorable statistic is an indicator function for the relative magnitude of each tower, and the subsequent statistics measure tilts and spreads within each tower.

We also examined the case where $g_0$ is a half-and-half mixture of a centered Gaussian with variance $\sigma^2 = 2$, and a $\delta$-spike at either $\mu = 0$ or $\mu = 2$. With a spike at $\mu = 0$, the most favorable statistics are not that different from the Hermite polynomials underlying \eqref{eq:poly}. However, the spike at $\mu = 2$ changes the picture considerably: now, the most stable statistic looks more like a hinge $T\p{\mu} = (\mu - 2)_-$.

Of course, it is not immediately clear how these non-Gaussian most favorable families are relevant to practical data analysis. A strength of the log-polynomial model \eqref{eq:poly} is that it is most favorable near \emph{any} Gaussian carrier $g_0$; meanwhile, the statistics presented in Figure \ref{fig:nongauss} depend on exact knowledge of $g_0$. However, whether or not they are immediately useful, the results presented in Figure \ref{fig:nongauss} can help us gain a better feeling for the flavor of density deconvolution, and for how the shape of the carrier $g_0$ affects the nature of the problem.

\section{Adaptive Local Minimaxity}
\label{sec:minimax}

In the previous section, we showed that there exist most favorable statistics $T$ that capture most of the relevant information in the $X$-sample in terms of our relative efficiency criterion. In this section, we show how to translate this result into a more classical minimax framework. Beyond being of independent interest, the minimax results provided here may also serve as additional evidence that the relative efficiency framework proposed above is ``natural'' in the sense that it helps motivate good statistical procedures.

In order to describe the asymptotics of Efron's method, we need a certain amount of formalism. To this end, we focus on the local perturbation model described in Section \ref{sec:setup}, and show that it induces a sequence of statistical problems that converge to a locally asymptotically normal experiment \citep{lecam1960locally}. This connection then lets us draw from the extensive literature on estimation in the Gaussian sequence model.

Here, we focus on the case where the carrier $g_0$ is Gaussian; this lets us cut down on linear algebra and to give closed-form bounds for the minimax shortfall of Efron's method. Following this choice, we use notation $g_\sigma$ instead of $g_0$ for the carrier; here, $\sigma^2$ denotes the variance of the $\mu_i$. We note, however, that exactly the same arguments can be used to derive the limiting risk of any most-favorable family of the type considered in Theorem \ref{theo:most_favorable}; the minimax performance of Efron's method then depends on the decay rate of the spectrum of $P_{g_0}$ \eqref{eq:lin}. We expand on this connection in Section \ref{sec:minimax_general}.

Following our discussion in Section \ref{sec:setup}, we are interested in a sequence of models
\begin{equation}
\label{eq:perturb}
\gn\p{\mu} = g_\sigma^M\p{\mu} \exp\sqb{\frac{1}{\sqrt{n}} \tau\p{\mu} - \psi_n},  \ \ \mu \in \Omega_M = [-M, \, M],
\end{equation}
where we eventually want to take $M \rightarrow \infty$.
We assume that the tilting function $\smash{\tau : \RR \rightarrow \RR}$ is a generic smooth function of $\mu$ contained within an ellipsoid defined with respect to its Hermite expansion
\begin{equation}
\label{eq:ellipse}
\tau\p{\mu} \in A^\sigma_{\kappa, \, C}, \ \ A^\sigma_{\kappa, \, C} = \cb{\sum_{j = 1}^\infty \kappa^{2j} \angles{\tau\p{\cdot}, \, H_j\p{\frac{\cdot}{\sigma}}}_{g_\sigma}^2 \leq C^2}.
\end{equation}
Here, the parameters $\kappa$ and $C$ let us tune the shape of the ellipsoid, and the inner product notation is short for
\begin{equation}
\label{eq:angles}
\angles{\tau\p{\cdot}, \, H_j\p{\frac{\cdot}{\sigma}}}_{g_\sigma} = \int_\RR \tau\p{\mu} H_j\p{\frac{\mu}{\sigma}} g_\sigma\p{\mu} d\mu.
\end{equation}
Our goal is to get a good estimate $\hgn$ for $\gn$ from $n$ independent samples $X_i$; we measure loss in terms of the re-scaled Kullback-Leibler divergence \eqref{eq:deviance}. Throughout our analysis, we will assume that $\hgn$ is chosen such as to ensure that this loss is finite.

Our key result is that, given an appropriate choice of $p$, Efron's polynomial log-density model \eqref{eq:poly} attains quasi-minimax performance simultaneously over this class of problems, regardless of our choice of $\sigma^2$ and $\kappa$.
We note that the constant on the right-hand side of \eqref{eq:minimax} can be quite small. For example, if we take the carrier to be standard Gaussian $\sigma^2 = 1$, and set the ellipse shape parameter to $\kappa = 2$, we have $\alpha_{\sigma, \, \kappa} \approx 1.7$ and $\beta_{\sigma, \, \kappa}  = 3$. Thus, the polynomial log-density model \eqref{eq:poly} gets to within a factor 1.8 of the minimax risk.

\begin{theo}
\label{theo:minimax}
Suppose that we are trying to solve the sequence of problems defined by \eqref{eq:perturb}, \eqref{eq:ellipse}, and \eqref{eq:deviance}. Let $\hgnp$ denote maximum likelihood estimator in the polynomial log-density model \eqref{eq:poly} with $p$ chosen as in \eqref{eq:pchoice}, and let $\hgns$ be the minimax optimal estimator over the regularity class defined in \eqref{eq:regularity_class}. Then, there is a constant $C_{\sigma, \, k}$ such that, for $C \geq C_{\sigma, \, \kappa}$,
\begin{equation}
\label{eq:minimax}
\lim_{M \rightarrow \infty} \limsup_{n \rightarrow \infty} \cb{ \frac{\sup_{\tau \in A^\sigma_{\kappa, \, C}}  {L_n\p{\hgnp}}} {\sup_{\tau \in A^\sigma_{\kappa, \, C}}  {L_n\p{\hgns}}} } \leq \frac{\beta_{\sigma, \, \kappa}}{\alpha_{\sigma, \, \kappa}},
\end{equation}
where the constants $\alpha_{\sigma, \, \kappa}$ and $\beta_{\sigma, \, \kappa}$ are given in \eqref{eq:alpha} and \eqref{eq:beta} respectively.
\end{theo}

To establish Theorem \ref{theo:minimax}, we first need to derive the worst-case risk of the minimax estimator $\hgns$. The key step in the proof is to show that the estimation problem outlined above converges to an elliptically constrained Gaussian sequence model. Now, as shown by \citet{pinsker1980optimal}, linear estimators get to within a constant factor of the minimax risk for this class of problems, and so the minimax estimation problem reduces to a linear estimation problem that can be solved directly; see \citet{johnstone2011gaussian} for a review. The factor $4/5$ in the constant \eqref{eq:alpha} is a bound obtained by \citet{donoho1990minimax} for the sub-optimality constant in Pinsker's result.

\begin{lemm}
\label{lemm:pinsker}
Under the conditions of Theorem \ref{theo:minimax}, for $C \geq C_{\sigma, \, \kappa}$,
\begin{equation}
\label{eq:pinsker}
\lim_{M \rightarrow \infty} \liminf_{n \rightarrow \infty} \cb{ \inf_{\hgns \in \Lambda_\sigma} \sup_{\tau \in A^\sigma_{\kappa, \, C}}  {L_n\p{\hgns}} } \geq \alpha_{\sigma, \, \kappa} \, C^{\frac{2\log\p{r_\sigma}}{\log{\p{r_\sigma} + \log\p{\kappa}}}},
\end{equation}
where $\Lambda^\sigma_\kappa$ is an $L_2$ regularity class defined in \eqref{eq:regularity_class}, and
\begin{equation}
\label{eq:alpha}
 \alpha_{\sigma, \, \kappa} = \frac{4}{5} \, \frac{r_\sigma^2 \p{\kappa - 1}}{\p{r^2_\sigma - 1}\p{r^2_\sigma \kappa - 1}} \p{\frac{\p{r_\sigma^2 \kappa - 1}\p{r^2_\sigma \kappa^2 - 1}}{r_\sigma^2 \kappa \p{\kappa - 1}}}^{\frac{2\log\p{r_\sigma}}{\log\p{r_\sigma} + \log\p{\kappa}}},
\end{equation}
and $r_\sigma$ is short-hand for the constant
\begin{equation}
\label{eq:r}
r_\sigma^2 = \p{1 + \sigma^2} \big/ {\sigma^2}.
\end{equation}
\end{lemm}

We can also derive the risk of $\hgnp$ using the machinery developed for proving Lemma \ref{lemm:pinsker}. Comparing \eqref{eq:pinsker} and \eqref{eq:poly_risk}, we see that both estimators have the same dependence on the signal scale $C$, and only differ by a constant function of $\sigma^2$ and $\kappa$.

\begin{lemm}
\label{lemm:poly_risk}
Suppose that the true density $\gn$ is as in the statement of Theorem \ref{theo:minimax}, and we estimate it using the exponential family model $\hgnp$ \eqref{eq:poly} with
\begin{equation}
\label{eq:pchoice}
p = \max\cb{2, \, \left\lceil \frac{\log\p{C / \sigma}}{\log\p{r_\sigma \kappa}} \right\rceil - 1}.
\end{equation}
Then, for $C \geq C_{\sigma, \, \kappa}$,
\begin{equation}
\label{eq:poly_risk}
\lim_{M \rightarrow \infty} \limsup_{n \rightarrow \infty} \cb{ \sup_{\tau \in A^\sigma_{\kappa, \, C}}  {L_n\p{\hgnp}} } \leq \beta_{\sigma, \, \kappa} C^{\frac{2\log\p{r_\sigma}}{\log\p{r_\sigma} + \log\p{\kappa}}},
\end{equation}
with
\begin{equation}
\label{eq:beta}
\beta_{\sigma, \, \kappa} = \p{1 + r^2_\sigma} \sigma^{\frac{2\log\p{\kappa}}{\log\p{r_\sigma} + \log\p{\kappa}}}.
\end{equation}
\end{lemm}

\subsection{Density Deconvolution and the Gaussian Sequence Model}
\label{sec:gauss_sequence}

Here, we outline the proof of Lemma \ref{lemm:pinsker} by showing how our density estimation problem converges to a Gaussian sequence model.
Throughout our analysis, we assume that our density estimate is in the class $\Lambda^\sigma_\kappa$ defined by the following relation
\begin{equation}
\label{eq:regularity_class}
\hg^*_{n} \in \Lambda^\sigma_\kappa, \ \ \Lambda^\sigma_\kappa = \cb{g : \sum_{j = 1}^\infty \kappa^{2j} \angles{\log\p{\frac{g\p{\cdot}}{g_\sigma^M\p{\cdot}}}, \, H_j\p{\frac{\cdot}{\sigma}}}_{g_\sigma^M} < L^2}
\end{equation}
for some large constant $L^2 > C^2$, along with the constraint that all integrals defined above be finite and well-defined.

Thanks to our regularity assumption, we can use notation from \eqref{eq:angles} to define ``Hermite coefficients''
\begin{align}
&\gamma_j = \angles{\tau\p{\cdot}, \, H_j\p{\frac{\cdot}{\sigma}}}_{g_\sigma}.
\end{align}
Meanwhile, given our sample $X_1, \, ..., \, X_n$, we can also define empirical Hermite coefficients as
\begin{equation}
\label{eq:hermite_coeff}
Z^{(n)}_j = \frac{1}{\sqrt{n}} \p{\frac{1 + \sigma^2}{\sigma^2}}^{j/2} \sum_{i = 1}^n H_j\p{\frac{X_i}{\sqrt{1 + \sigma^2}}}.
\end{equation}
Because $\smash{\gn}$ converges to the carrier $g_0$ as $n$ gets large, we can show that for any finite set of indices $\smash{\cb{j_i}}$, the $\smash{Z^{(n)}_{j_i}}$ are asymptotically jointly normal with
\begin{align}
\label{eq:gauss_cov}
&\limn \EE{Z^{(n)}_{j_i}} = \gamma_{j_i} + o_M(1), \\
&\limn \Cov{Z^{(n)}_{j_{i}}, \, Z^{(n)}_{j_{i'}}} = \delta_{\cb{j_i = j_{i'}}} r_\sigma^{2{j_i}} + o_M(1),
\end{align}
where we again used the notation $\smash{r_\sigma^2 = \p{1 + \sigma^2}/{\sigma^2}}$.

These observations suggest that deriving a good estimator $\smash{\hgn}$ for the density $\smash{\gn}$ is related to finding the mean of the Gaussian sequence with covariance structure \eqref{eq:gauss_cov}. The following lemma makes this connection explicit, thus reducing our setting to a well-understood problem, namely estimating a Gaussian sequence model with elliptical constraints under squared-error loss.

\begin{lemm}
\label{lemm:gauss_sequence}
For large $M$, the limiting minimax risk for the sequence of density deconvolution problems defined in Lemma \ref{lemm:pinsker} converges to the minimax risk for estimating the mean of a Gaussian sequence $Z_j$ for $j = 1, \, 2, \, ...$ with
\begin{equation}
\label{eq:gauss_sequence}
\EE{Z_j} = \gamma_j, \ \Cov{Z_j, \, Z_k} = \delta_{\cb{j = k}} r_\sigma^{2j}
\end{equation}
under the loss
\begin{equation}
\label{eq:loss_gauss}
L\p{\hgamma\p{Z}} = \sum_{j = 1}^\infty \p{\hgamma\p{Z} - \gamma}^2,
\end{equation}
and the constraint $\smash{\gamma \in  \ell_2^\kappa\p{C} := \cb{\gamma : \sum_{j = 1}^\infty \kappa^{2j} \gamma_j^2 \leq C^2}}$. In other words,
\begin{align}
&\lim_{M \rightarrow \infty} \limn \cb{\inf_{\hg \in \Lambda^\sigma_\kappa} \sup_{\tau \in A^\sigma_{\kappa, \, C}} L_n \p{\hg^{(n)} \p{Z^{(n)}}}} \\
\notag
&\ \ \ \ \ \ = \inf_{\hgamma \in  \ell_2^\kappa\p{C}} \sup_{\gamma \in \ell_2^\kappa\p{C}} L\p{\hgamma\p{Z}}.
\end{align}
\end{lemm}

Now, this class of Gaussian sequence models can be well estimated using linear rules: \citet{pinsker1980optimal} established that the risk of the best linear rule is within a constant factor of the minimax risk; furthermore, \citet{donoho1990minimax} showed that this constant is less than 5/4. Thus, it suffices to find the risk of the best linear rule of the form $\hgamma^{L}_j = c_j Z_j$ for some constants $c_j$. Now, in the Gaussian limit, the worst-case risk of the minimax linear rule $\hgamma^L$ has a simple form \citep[e.g.,][ Chapter 5.1]{johnstone2011gaussian}:
\begin{equation}
\label{eq:rl}
R^L :=  \sup_{\gamma \in \ell_2^\kappa\p{C}} \EE{L\p{\hgamma^L}} = \sum_{j = 1}^\infty  r_\sigma^{2j} \p{1 - \frac{\kappa^j}{\mu_C}}_+,
\end{equation}
where $\mu_C$ is implicitly defined by the relation
\begin{equation}
\label{eq:implicit}
\sum_{j = 1}^\infty r_\sigma^{2j} \kappa^j \p{\mu_C - \kappa^j}_+ = C^2.
\end{equation}
Since we know that the minimax risk of $R^*$ is bounded from below by $4 R^L / 5$,  the proof of Lemma \ref{lemm:pinsker} reduces to algebra; the remaining steps are carried out in the appendix.

\subsection{Extension to General Carriers}
\label{sec:minimax_general}

In an effort to reduce notational burden, the above argument focused on the Gaussian carrier case, i.e. $g_0 = g_\sigma$ for some $\sigma > 0$. In this section, we briefly outline the technical ideas needed to prove analogous local minimaxity results in the neighborhood of general carriers $g_0$.
The proof of Theorem \ref{theo:minimax} relied on showing that large-sample statistical inference in the local density deconvolution model reduces to a study of the asymptotically normal ``empirical Hermite coefficients'' $\smash{Z_j^{(n)}}$ defined in \eqref{eq:hermite_coeff}. A similar result also holds in the case of general $g_0$; however, the statistics $\smash{Z_j^{(n)}}$ now take on a more general form.

Let $\smash{\cb{T_j(\mu)}_{j = 1}^\infty}$ be the statistics comprising the most-favorable family around $g_0$, and denote the relative efficiency of the $j$-dimensional most-favorable family by $\rho_j$. Then, inference can be asymptotically framed in terms of the moments 
\begin{equation}
\label{eq:general_coeff}
Z_j^{(n)} = \frac{1}{\sqrt{n}} \, \rho_j^{-1} \, \sum_{i = 1}^n U(X_i), \ \with \ U(x) := \frac{\p{\phi * \p{T_j \, g_0}}\p{x}}{f_0\p{x}}.
\end{equation}
The following lemma shows that the $Z_j^{(n)}$ in fact have the desired limiting distribution. An analogue to Lemma \ref{lemm:gauss_sequence} then follows directly.

\begin{lemm}
\label{lemm:general_coeff}
Suppose that the conditions of Theorem \ref{theo:most_favorable}, and that the perturbation function $\tau$ satisfies $\smash{\int_\Omega \tau^2\p{\mu} g_0\p{\mu} \ d\mu < \infty}$. Then, then the $Z_j^{(n)}$ defined in \eqref{eq:general_coeff} are asymptotically normal over any finite set of indices and, for any $j, \, j' \in \NN$,
\begin{align}
&\limn \EE{Z_j^{(n)}} = \int_\Omega  T_j\p{\mu} \, \tau\p{\mu} \, g_0\p{\mu} \ d\mu, \\
&\limn \Cov{Z_j^{(n)}, \, Z_{j'}^{(n)}} = \delta_{\cb{j = j'}} \, \rho_j^{-1}.
\end{align}
\end{lemm}

Given this result, we can use Pinsker's theorem to compute tight bounds for the minimax risk of density deconvolution near $g_0$  as a function of the most-favorable relative efficiency coefficients $\rho_j$ (i.e., the spectrum of $P_{g_0}$ defined in \eqref{eq:lin}). We note that Lemma \ref{lemm:general_coeff}---and in fact the whole machinery of using a Gaussian sequence model to understand the behavior of maximum likelihood estimation in the most-favorable family---only depends on the square-integrability condition $\smash{\int_\Omega \tau^2\p{\mu} g_0\p{\mu} \ d\mu < \infty}$ and on the compactness of $P_{g_0}$. However, in order for the resulting minimaxity properties to be any good, we need for the spectrum of $P_{g_0}$ to decay fast.

\section{Two Optimality Theories for Density Deconvolution}
\label{sec:examples}

The main contribution of our paper is a local optimality theory for density deconvolution that helps us understand and justify Efron's method, i.e., maximum likelihood estimation in a model of the form \eqref{eq:expfam}. This line of work is in contrast to the classical optimality theory based on kernel estimators. Notable contributions include the pioneering work of \citet{carroll1988optimal} and \citet{stefanski1990deconvoluting}, the analysis of \citet{fan1991optimal} that elucidates the connection between the decay rate of the Gaussian characteristic function and the difficulty of density deconvolution, the strong quasi-minimaxity results of  \citet{efromovich1997density}, as well as several other papers cited in the introduction.

A key difference between these theories is that kernel-based methods are quasi-optimal in terms of an integrated squared error (ISE) loss criterion whereas our theory is framed in terms of the Kullback-Leibler (KL) loss, defined respectively as
\begin{align}
\label{eq:ISE}
&\text{ISE}\p{g, \, \hg} = \int \p{\hg\p{\mu} - g\p{\mu}}^2  d\mu, \\
\label{eq:KL}
&\text{KL}\p{g, \, \hg} = \int \log\p{\frac{g\p{\mu}}{\hg\p{\mu}}} g\p{\mu} d\mu.
\end{align}
Moreover, our theory aims to detect \emph{local} perturbations of $g_0$, whereas the ISE criterion is usually applied \emph{globally}.

From a scientific point of view, the value of an optimality theory depends on the relevance of the induced estimators to answering real-world questions. Kernel-based methods have a good track record for solving some classic problems, such as in-season baseball prediction problem introduced by \citet{efron1975data}; see, e.g., \citet{brown2008season}. In fact, these methods have been shown to approach the Bayes risk for estimating the posterior mean $\EE{\mu \cond X = x}$ \citep{brown2009nonparametric}.

In this section, however, we present examples of natural scientific questions for which Efron's method provides substantially better answers than kernel methods. We hope that these examples will convince the reader that a KL-based optimality theory for density deconvolution is, if nothing else, worthy of further study. For our experiments, we used the \texttt{R}-package \texttt{decon} \citep{wang2011deconvolution} for kernel-based estimation.

\subsection{Detecting the Fraction of Weakly Associated Genes in a Microarray Study}

Our first example simulates a classic biological application: gene expression profiling using microarrays. At a high level, the goal is to estimate the difference in the expression levels of different genes for two distinct cell populations (e.g., breast cancer cells vs. healthy cells). After some pre-processing, each gene can be associated with a test statistic that has a standard normal distribution under the null hypothesis that the gene's expression levels do not differ accross the two groups; the resulting statistical problem is described in detail by \citet{efron2001empirical} and \citet{tusher2001significance}.

Here, we follow the \emph{structural model} analysis of \citet{efron2004large} who showed that, to reasonable approximation, we can model the gene-wise test statistics $X_i$ as
$ X_i \sim \nn\p{\mu_i, \, 1}, $
where $\mu_i$ describes the true association between the $i$-th gene and the condition of interest. The $i$-th ``exact'' null hypothesis is that $\mu_i = 0$. However, as argued in, e.g., Chapter 6 of \citet{efron2010large}, this exact null may not always be scientifically relevant. It seems likely that most genes have a small but non-zero true association $\mu_i$; the goal of the statistician is then not to identify the non-zero $\mu_i$ but rather to identify the \emph{large} associations $\mu_i$.

\begin{figure}[t]
\centering
\includegraphics[width = 0.6\textwidth]{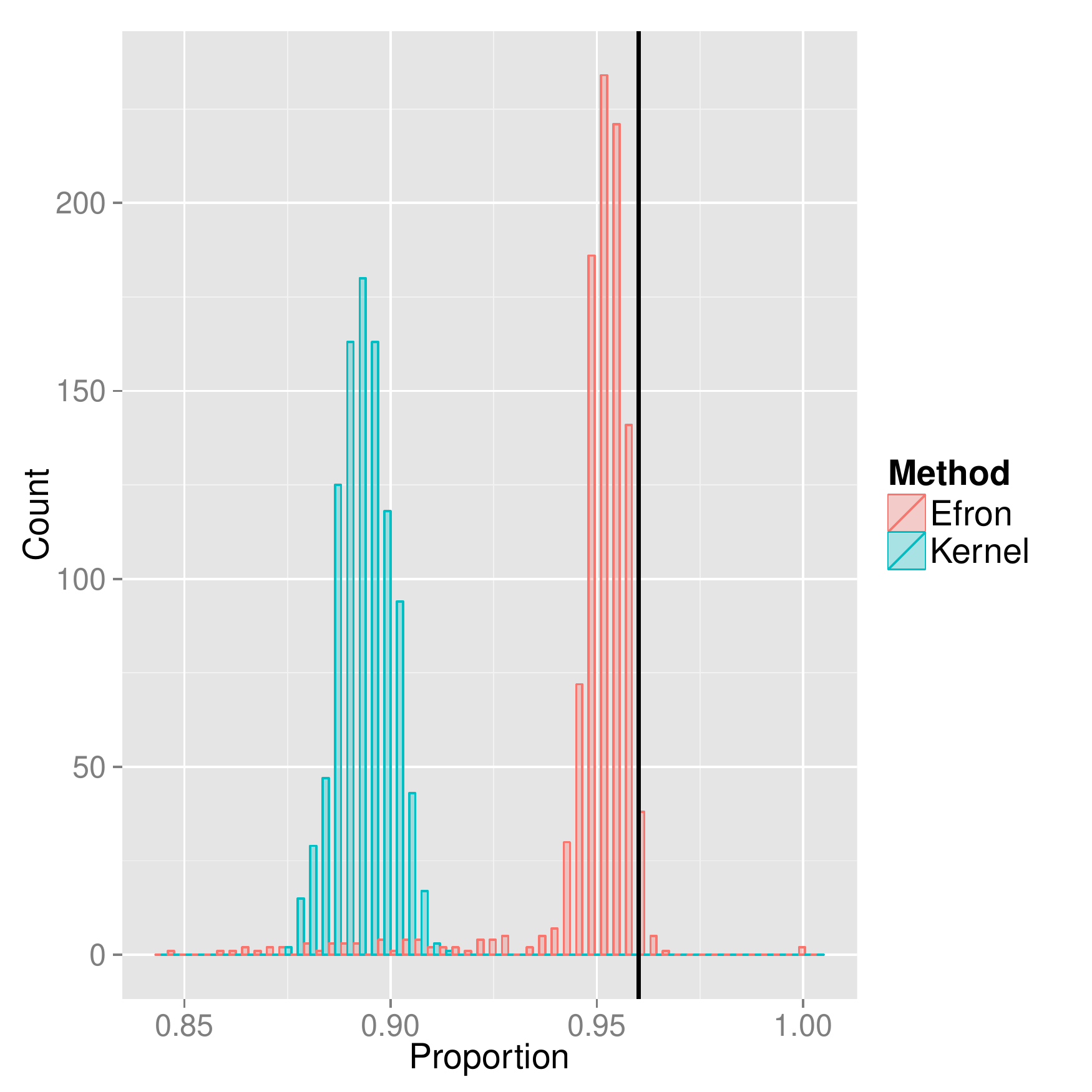}
\caption{Gene expression profiling simulation. The goal is to estimate $\int_{-2}^2 g\p{\mu} \, d\mu$; here, the correct answer (0.96) is indicated by a thick vertical line. Over 1,000 simulation replicates, Efron's method usually performs substantially better than the kernel-based alternative.}
\label{fig:gene}
\end{figure}

Motivated by this setup, suppose that---given the power afforded by our sample size---we decide that having $\mu_i \in [-2, \, 2]$ qualifies as a  ``small'' association, and we want to estimate the fraction of genes whose association is in this range. In our framework, assuming that the gene associations $\mu_i$ are drawn from a distribution with density $g(\cdot)$, we want to estimate $\int_{-2}^2 g\p{\mu} \, d\mu$. In Figure \ref{fig:gene}, we ran 1,000 replicates of such a simulation with $n = 5,000$ genes each, where the true density $g$ was given by
$$ g\p{\mu} = 0.95 \cdot \frac{1}{4} \p{2 - \abs{\mu}}_+ + 0.05 \cdot \frac{1}{20} 1\p{\cb{\abs{\mu} \leq 10}}. $$
We then generated $X$-statistics from the structural model $X \sim \nn\p{\mu, \, 1}$, and generated estimates $\hg$ both using the kernel method \texttt{decon} and maximum likelihood in the family \eqref{eq:poly} with $p = 4$.

In this simulation, the correct answer was $\int_{-2}^2 g\p{\mu} \, d\mu = 0.96$. Efron's method \eqref{eq:poly} does substantially better than the kernel method, with the exception of a few cases where it vastly underestimates the bulk of $g(\cdot)$. In practice, \citet{efron2014bayes} recommends the use of regularization to stabilize the estimators and mitigate such finite sample effects. Here, however, we did not use any regularization in order to keep our experiments as simple as possible.

\subsection{Identifying Safe Neighborhoods}

Our second example is based on a real dataset: the communities and crime unnormalized dataset from the UCI Machine Learning Repository \citep{UCI}, which tallies crime data from multiple US communities in 1995; we restricted our analysis to the $n = 1,192$ communities with population greater than 20,000 and non-missing crime data.

Our goal was to understand the number of non-violent crimes per 100 people. In our sample of communities, the average non-violent crime rate was around 5 per 100 inhabitants, while 6\% of communities achieved a non-violent crime rate of 2\%. We define a community with a rate of less than 2\% as \emph{safe}. The dataset was large enough that we could accurately estimate per-community crime rates; to test our deconvolution methods, we made the problem harder by down-sampling the data and then seeking to recover the correct answer.

\begin{figure}[t]
\centering
\includegraphics[width = 0.7\textwidth]{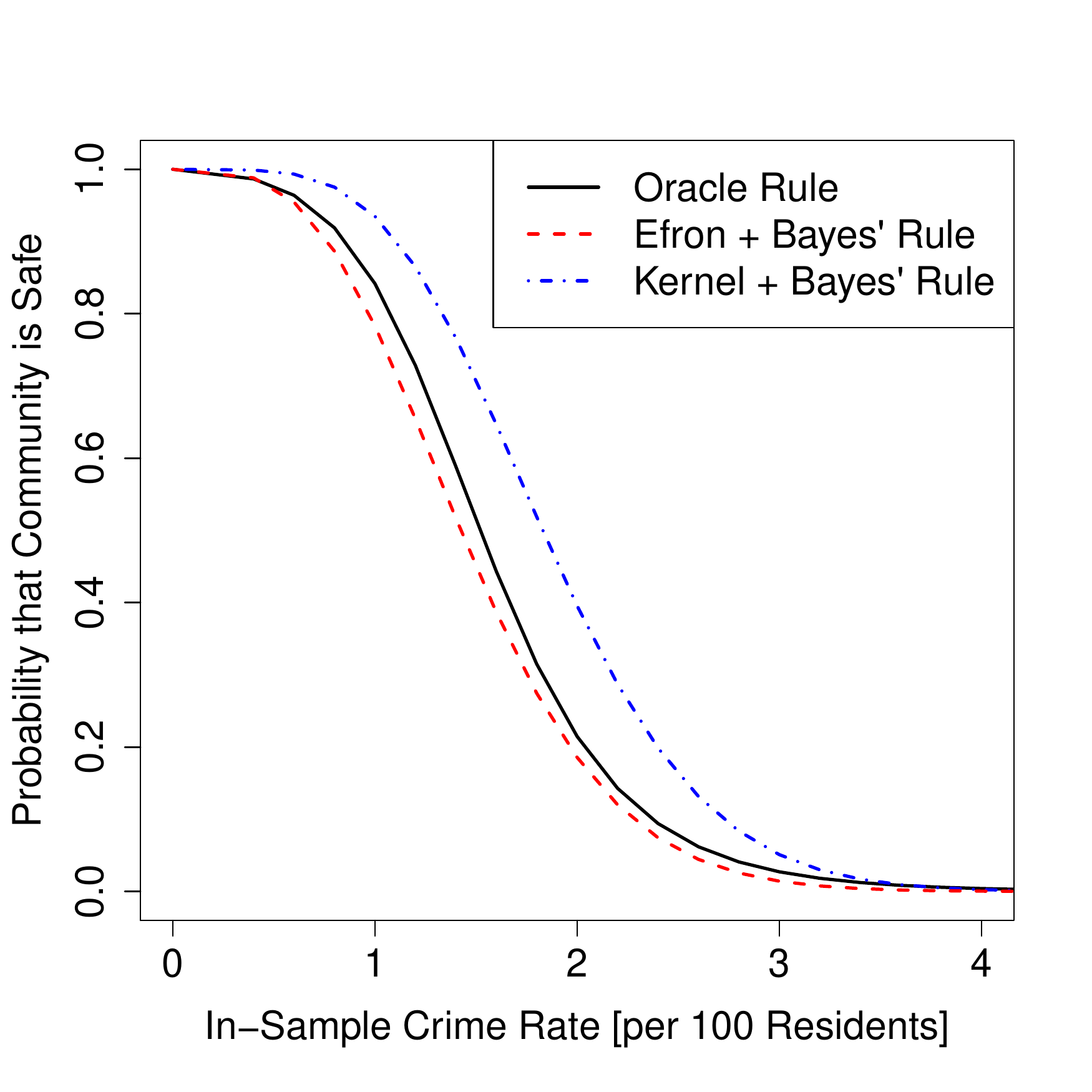}
\caption{Crime prediction example. The goal is to predict the probability that a community is safe ($p_i \leq 0.02$) given a crime-rate estimate $\hp_i$ obtained by interviewing $B = 500$ randomly selected people. Efron's method provides a closer approximation of the oracle rule than the kernel method.}
\label{fig:crime}
\end{figure}

To down-sample the data in a mathematically principled way, we made the practically somewhat implausible assumption that each person was victim of at most 1 crime in 1995. Then, we can imagine assembling a crime dataset by interviewing $B = 500$ randomly selected people per community, and counting the number $N_i$ of interviewees in the $i$-th community that have been victims of a non-violent crime. We can easily simulate the outcome of such an interview using the available data by hypergeometric sampling:
$$ N_i \sim \text{Hyper}\p{B, \, \text{Crimes in community $i$}, \, \text{Population $i$}}. $$
Our statistical task is to estimate which communities are \emph{safe} (i.e., have a rate of less than 2\%) based on statistics $N_i$ collected in different communities.
Let $p_i$ denote the true crime rate in the $i$-th community, and let $\hp_i = N_i / B$. By a standard variance stabilizing argument, we can verify that
$$ \sqrt{\hp_i} \overset{\cdot}{\sim} \nn\p{\sqrt{p_i}, \, \frac{1}{B}}. $$
Thus, to estimate $\PP{p_i \leq 0.02 | \hp_i}$, we can first apply a square-root transform to the data, then use any method to estimate the density $g$ of the $\sqrt{p_i}$, and finally apply Bayes' rule to get our desired quantity.

Figure \ref{fig:crime} shows results for the kernel method and Efron's method; both density deconvolution methods were implemented exactly the same way as for the gene expression example. The oracle rule was produced by fitting a spline logistic regression of $1\p{\cb{p_i \leq 0.02}}$ against $\hp_i$. Efron's method again substantially outerforms the alternative. For example, if we ask about the probability that a community is safe given that $\hp_i = 0.02$, Efron's method tells us that this probability is $18.6\%$ whereas the kernel method tells us $39.6\%$. The correct answer, produced by the oracle fit, is $21.5\%$.

\subsection{Fitting the Mode vs. Fitting the Support}

\begin{figure}[t]
\centering
\includegraphics[width = 0.7\textwidth]{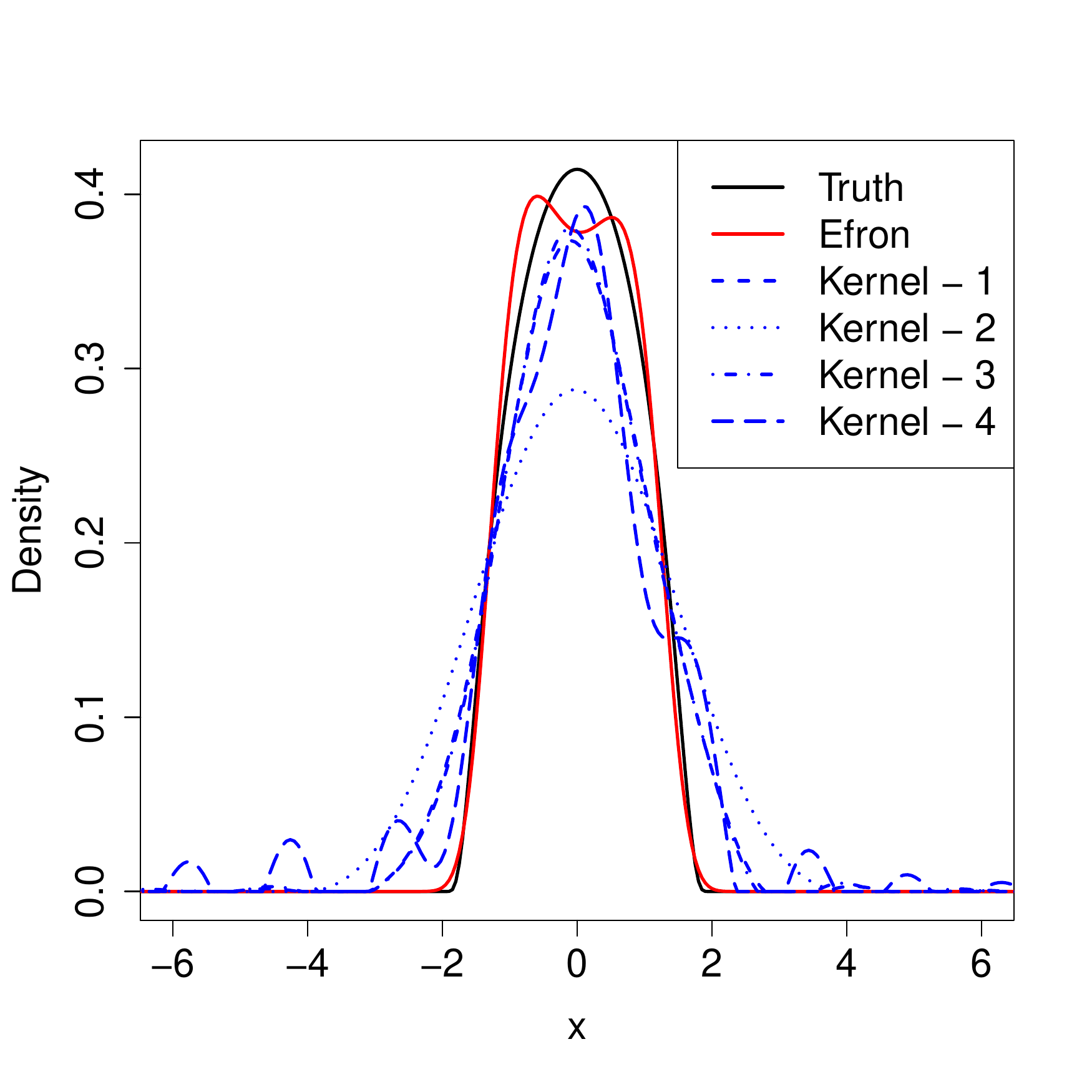}
\caption{A comparison of Efron's method with a log-polynomial model \eqref{eq:poly} with $p = 4$ with kernel-based density deconvolution as implemented in \texttt{decon}. We tried four different bandwidth choices for the kernel method.}
\label{fig:cmp}
\end{figure}

To provide some deeper insight into the differences between the two density deconvolution methods under consideration, we end our experimental section with an in-depth look at a simple density estimation problem with a large sample size: we drew $n = 100,000$ observations from the generative model \eqref{eq:setup} with a density $g$ of the form
$$ g(\mu) \propto \exp\sqb{- 1 /\p{1 + \frac{4}{\mu^2}}} \text{ for } -2 < \mu < 2, \text{ and } g(\mu) = 0 \text{ else.} $$
We note that this density is not contained in the span of either \eqref{eq:poly} or the basis functions implicitly used by \texttt{decon}. 
Results are shown in Figure \ref{fig:cmp}.

To give the kernel density method a good chance of performing well, we tried four different bandwidth-selection algorithms provided by \texttt{decon}: {\bf (1)} a closed form approximation to the boostrap recommended by \citet{delaigle2004bootstrap} (this method is default, and was used for the other experiments), {\bf (2)} a rule of thumb by \citet{fan1991optimal}, {\bf (3)} an approximate ISE minimzer by \citet{stefanski1990deconvoluting}, and finally {\bf (4)} a hand-selected bandwidth to demonstrate the effect of mild under-smoothing. For Efron's method we, as usual, set $p = 4$ in \eqref{eq:poly}.

We notice that all methods do a roughly equivalent job of estimating the target density $g$ near its mode. However, Efron's method is much more accurate in estimating the tails and support of $g$. This observation is not entirely surprising in light of the loss functions that were used to motivate each method. The ISE loss \eqref{eq:ISE} only requires $\hg$ to be reasonable close to $g$ over the domain of $\mu$, but does not impose particularly harsh penalties for oscillatory behavior in the tails. In contrast, the KL loss places more attention on getting the tails and support of the distribution right. It thus appears that kernel methods built on ISE-optimality theory are unreliable for answering scientific questions that depend on understanding the tail-behavior of $g$, whereas methods based on KL-optimality may perform better.

\section{Discussion}
\label{sec:discussion}

The work of \citet{efron2014bayes,efron2014two} presents a surprising challenge to the theory of density deconvolution. From, e.g., \citet{efromovich1997density}, it may appear that the problem of density deconvolution with Gaussian errors is completely solved---and that kernel density estimators are optimal for the task. And yet, \citet{efron2014two} found that his maximum likelihood method vastly out-performed methods related to kernel density estimators for several realistic problems (Efron calls this latter approach ``$f$-estimation'').

The goal of this paper was to understand why Efron's method could work better than kernel density estimators despite the theoretical guarantees available for the latter. To do so, we introduced a perturbation model inspired by Le Cam's local asymptotic normality theory, and showed that Efron's method is quasi-optimal in this setup for deviance loss. Under the assumption that deviance loss comes closer to describing the ``real loss function'' of a practitioner than the integrated squared error loss used to establish optimality properties of kernel density estimators, our results can be seen as moving towards an explanation for the empirical success of Efron's method.

More broadly, our results highlight a surprisingly strong connection between non-parametric density deconvolution and classical likelihood theory. We found that---at least locally---there are some directions along which we can accurately estimate a signal and maximum-likelihood estimation is quasi-optimal for this task; meanwhile, there are other directions in which estimation is hopeless and the minimax strategy is to ignore them. From a practical perspective, this connection appears rather reassuring, as low-dimensional large-sample maximum-likelihood estimation has proven to be one of the most consistently successful ideas in applied statistics.

\setlength{\bibsep}{0.5pt plus 0.3ex}
\bibliographystyle{plainnat}
\bibliography{references}


\begin{appendix}
\section{Proofs}
\label{sec:proofs}

\subsection*{Proof of Lemma \ref{lemm:relative}}

It suffices to prove the second conclusion, as the alternative expression for the relative efficiency coefficient given in \eqref{eq:rel_mult2} can directly be verified to be transformation invariant. We already know that $\ii_\mu\p{\eta} = \Var[\eta]{t(\mu)}$, and so all we really need to do is to check that
$$ \ii_X\p{\eta} = \Var[\eta]{\EE[\eta]{t(\mu) \cond X}}. $$
Although we could also get to this answer by drawing analogies to the work of \citet{louis1982finding}, we will take a direct approach here. In the univariate family with statistic $t$, the expected score function  is
\begin{align*}
\frac{d}{d\eta} \log f_\eta \p{x}
&= \frac{d}{d\eta} \log \int K\p{\mu, \, x} g_\eta\p{\mu} d\mu \\
&= \frac{\int K\p{\mu, \, x} g_\eta'\p{\mu} d\mu} {\int K\p{\mu, \, x} g_\eta\p{\mu} d\mu} \\
&= \frac{\int K\p{\mu, \, x} g_\eta\p{\mu} \p{t(\mu) - \EE[\eta]{t(\mu)}}d\mu} {\int K\p{\mu, \, x} g_\eta\p{\mu} d\mu} \\
&= \EE[\eta]{t(\mu) \cond X = x} - \EE[\eta]{t\p{\mu}}.
\end{align*}
Thus, we conclude that
\begin{align*}
\ii_X\p{\eta}
&= \Var[\eta]{\frac{d}{d\eta} \log f_\eta \p{x}} \\
&= \Var[\eta]{\EE[\eta]{t(\mu) \cond X = x}}.
\end{align*}

\subsection*{Proof of Lemma \ref{lemm:mult}}

By Bayes' rule, we know that the conditional density of $\mu$ given $x$ is
$$ g\p{t\p{\mu} \cond X = x} = \frac{K\p{x, \, \mu} g\p{\mu}}{f\p{x}}; $$
from this, it directly follows that
$$ \Var{\EE{t\p{\mu} \cond X}} = \int_\Omega  t^2\p{\mu}  K^2\p{x, \, \mu} g^2\p{\mu}f^{-1}\p{x}  dx \, d\mu. $$ 
Similarly, we can check that
$$ \Var{t\p{\mu}} =  \int_\Omega  t^2\p{\mu}g\p{\mu} d\mu. $$
For a statistic $t\p{\mu}$ of the form $t = a \cdot T$, we can write $t^2 = a^\top (T T^\top) a$; with this notation, we recover the first result of Lemma \ref{lemm:mult}. The explicit solution given in \eqref{eq:prag1} and \eqref{eq:prag2} for $a^*$ in terms of the spectrum of $Q_T^\top \, M_T Q_T$ is standard.

\subsection*{Proof of Theorem \ref{theo:most_favorable}}

Because the relative efficiency coefficient is transformation invariant, we can without loss of generality pick a statistic $T$ for which
\begin{align}
\label{eq:constraint1}
&\int_\Omega T\p{\mu}  g_0\p{\mu} d\mu = 0, \\
\label{eq:constraint2}
&\int_\Omega T\p{\mu} T^\top\!\!\p{\mu} g_0\p{\mu} d\mu = I_{p \times p}. 
\end{align}
Given such a choice, the relative efficiency formula from Lemma \ref{lemm:mult} simplifies to
\begin{align}
\label{eq:objective}
\rho_0\p{T}
&= \lambda_{\min}\p{ \int_\Omega T\p{\mu} T^\top\!\!\p{\mu} K^2\p{x, \, \mu} g_0^2\p{\mu}f^{-1}\p{x}  dx \, d\mu}, \\
\notag
&= \lambda_{\min}\p{ \int_\Omega \p{\sqrt{g_0}\p{\mu} T\p{\mu}}  P_{g_0}\p{\mu, \, \mu} \p{\sqrt{g_0}\p{\mu}T}^\top d\mu}, 
\end{align}
where $\lambda_{\min}(A)$ denotes the smallest eigenvalue of a linear operator $A$. Minimizing the objective \eqref{eq:objective} subject to the constraint \eqref{eq:constraint2} is a standard spectral analysis problem. By the Courant-Fischer-Weyl maximin theorem as stated in, e.g., \citet{shawe2005eigenspectrum}, we find that because $P_{g_0}$ is both self-adjoint and compact (and thus also completely continuous),
$$ \max\cb{\rho_0\p{S} : S \in L_2\p{\Omega}^p, \, \int_\Omega S\p{\mu} S^\top\!\!\p{\mu} g_0\p{\mu} d\mu = I_{p \times p}} = \lambda_p, $$
where $\lambda_p$ is the $p$-the eigenvalue of $P_{g_0}$; moreover, this maximum is attained by setting $S_j\p{\mu} = \zeta_j\p{\mu} / \sqrt{g_0(\mu)}$ for $j = 1, \, ..., \, p$, where the $\zeta_j$ are the leading eigenvectors of $P_{g_0}$. We note that $\zeta_j\p{\mu} \in L_2\p{\Omega}$ because $P_{g_0}$ is compact; $S_j$ is then also in $L_2\p{\Omega}$ because $g_0$ is bounded away from 0 on $\Omega$. Finally, because $g_0$ and $K$ are continuous, $P_{g_0} S_j$ is continuous and so $S_j$ must also be continuous.

Now, we still need to deal with the constraint \eqref{eq:constraint1}. Thankfully, we can verify that all the eigenvalues of $P_{g_0}$ are bounded by 1, and that $\zeta_1\p{\mu} := \sqrt{g_0\p{\mu}}$ is an eigenfunction of $P_{g_0}$ with eigenvalue 1. By orthogonality of the spectrum of $P_{g_0}$, we then see that
$$ \int_\Omega \frac{\zeta_j\p{\mu}}{\sqrt{g_0\p{\mu}}} g_0\p{\mu} d\mu = \int_\Omega \zeta_j\p{\mu} \zeta_1\p{\mu} = 0 $$
for all $j > 1$. Thus, the minimizer of \eqref{eq:objective} with both constraints \eqref{eq:constraint1} and \eqref{eq:constraint2} is given by $T_j\p{\mu} = \zeta_{j + 1}\p{\mu} / \sqrt{g_0(\mu)}$ for $j = 1, \, ..., \, p$, and the objective value \eqref{eq:objective} is the $(p+1)$-st eigenvalue of $P_{g_0}$. Moreover, if the spectrum of $P_{g_0}$ does not have repeated eigenvalues, the span of $T_1, \, ..., \, T_p$ maximizing our objective is unique. We note that, because $g_0$ and $S_j$ are continuous and $\Omega$ is compact, $\EE{g_0\p{\mu} \exp\p{\eta \cdot T\p{\mu}}}$ is finite for $\eta$ in a neighborhood of 0 and so our most-stable family is in fact well-defined.

\subsection*{Proof of Corollary \ref{coro:approx}}

Because $P_{g_0}$ is compact, we know by the eigenfunctions $\smash{\cb{\zeta_j}_{j = 1}^\infty}$ form a complete orthonormal basis for $L_2\p{\Omega}$; thus, we see that the statistics $\smash{\cb{T_j}_{j = 0}^\infty}$ also form a complete orthonormal basis for $L_2\p{\Omega}$ with inner product weighted by $g_0$ and can write\footnote{Without loss of generality, $\gamma_0 = 0$ since this term is absorbed by the normalization.}
$$ \tau\p{\mu} = \sum_{j = 1}^\infty \gamma_j \, T_j\p{\mu}, \ \ \int_\Omega \tau^2\p{\mu} g_0\p{\mu} = \sum_{j = 1}^\infty \gamma_j^2. $$
Given this notation, we set $\smash{\tau^{(p, \, *)} = \sum_{j = 1}^p \gamma_j \, T_j\p{\mu}}$.
Our goal is to show that, for any integer $J > p$, our choice of $\smash{\tau^{(p, \, *)}}$ satisfies the conclusion of Corollary \ref{coro:approx} under the assumption that $\gamma_j = 0$ for all $j > J$. Because this bound holds uniformly in $J$, we conclude that it also holds in the non-parametric case $\smash{\int_\Omega \tau^2\p{\mu} g_0\p{\mu} \leq C^2}$.

Following the above discussion, we now assume that $\smash{\tau\p{\mu} = \sum_{j = 1}^J \gamma_j \, T_j\p{\mu}}$, and write $\gamma$ for the parameter vector inducing $\tau$. The target loss is then
\begin{align*}
D_{KL} &\p{f^{(n)}_\tau, \, f^{(n)}_{\tau^{(p, \, *)}}} \\
&= \int_\Omega f_{\tau/\sqrt{n}}\p{x} \, \log\p{\frac{f_{\tau/\sqrt{n}}\p{x}}{f_{\tau^{(p, \, *)}/\sqrt{n}}\p{x}}} \, dx \\
&= \frac{1}{2n} \int_\Omega   \p{\gamma^{(p, \, *)} - \gamma}^\top \nabla^2 \log\p{ f_{\tau/\sqrt{n}}\p{x}} \p{\gamma^{(p, \, *)} - \gamma} \, f_{\tau/\sqrt{n}}\p{x}  \, dx \\
&\ \ \ \ \ \ \ \ \ \ \ \ \ \  + o\p{n^{-1}},
\end{align*}
where the derivative $\nabla^2$ is taken with respect to the parameters $\gamma$. We note that the first-order term depending on $\nabla \log(\cdot)$ integrates out to 0. Now, taking limits, we find that
\begin{align*}
\limn & n \, D_{KL} \p{f^{(n)}_\tau, \, f^{(n)}_{\tau^{(p, \, *)}}} \\
&= \frac{1}{2} \p{\gamma^{(p, \, *)} - \gamma}^\top
\int_\Omega \nabla^2 \log\p{f_0\p{\mu}} f_0\p{x} \, dx \
\p{\gamma^{(p, \, *)} - \gamma} \\
&= \frac{1}{2} \p{\gamma^{(p, \, *)} - \gamma}^\top \ii_X^{\Gamma_J\p{g_0}} \p{\gamma^{(p, \, *)} - \gamma},
\end{align*}
where $\smash{\ii_X^{\Gamma_J\p{g_0}}}$ denotes the Fisher information for estimating $\gamma$ from $X$-samples in the $J$-dimensional most-favorable family.
But now, $\smash{\Gamma_J\p{g_0}}$ is scaled such that $\smash{\ii_\mu^{\Gamma_J\p{g_0}}} = I_{J \times J}$. Thus, the limiting approximation error is equal to 
$$ \limn n \, D_{KL} \p{f^{(n)}_\tau, \, f^{(n)}_{\tau^{(p, \, *)}}} = \frac{1}{2} \Norm{\gamma^{(p, \, *)} - \gamma}_2^2 \, \rho\p{\p{\gamma^{(p, \, *)} - \gamma} \cdot \Gamma_j\p{g_0}}, $$
which, by Theorem \ref{theo:most_favorable}, can be bounded above by $\frac{1}{2} \,  C^2 \lambda_{p+2}\p{P_{g_0}}$.

\subsection*{Proof of Theorem \ref{theo:gauss_hard}}

Given our assumptions, we know that the carrier $g$ and the marginal density of the observations $f$ are given by
$$ g\p{\mu} = \frac{1}{\sigma} \, \varphi\p{\frac{\mu}{\sigma}}, \ \ f\p{x} = \frac{1}{\sqrt{1 + \sigma^2}} \, \varphi\p{\frac{x}{\sqrt{1 + \sigma^2}}}, $$
where both densities loop around if needed to accommodate the bounded domain.
For our proof, we begin by verifying that the
$$ \nu_j\p{\mu} = \sqrt{ \frac{1}{\sigma} \, \varphi\p{\frac{\mu}{\sigma}}} \,  H_j\p{\frac{\mu}{\sigma}} $$
are eigenfunctions of $P_{g}$ in the limit $M = \infty$. Then, for large but finite $M$, the $\nu_j$ are nearly eigenfunctions of $P_g$; meanwhile, the conditions of Theorem \ref{theo:most_favorable} are satistified and so our desired conclusion follows.

Now, to verify that the $\nu_j$ are eigenfunctions with $M = \infty$, we first note that $P_g$ is a compact kernel and so it does in fact admit a spectral decomposition, and second that
\begin{align*}
\int &\nu_j\p{\mu_1} P_{g} \p{\mu_1, \, \mu_2} \nu_k\p{\mu_2} d\mu_1 d\mu_2 \\
&= \int \p{\varphi * \p{\frac{1}{\sigma} \, \varphi\p{\frac{\cdot}{\sigma}} \cdot H_j\p{\frac{\cdot}{\sigma}}}}\p{x} \\
&\ \ \ \ \ \ \ \ \ \ 
 \p{\varphi * \p{\frac{1}{\sigma} \, \varphi\p{\frac{\cdot}{\sigma}} \cdot H_k\p{\frac{\cdot}{\sigma}}}}\p{x}
 \frac{dx}{f(x)}
\end{align*}
Now, focusing on the inner terms, we can check that
\begin{align*}
\frac{1}{\sigma} \, \varphi\p{\frac{\mu}{\sigma}} \cdot H_j\p{\frac{\mu}{\sigma}}
= \frac{1}{\sqrt{j!}}  \, \frac{1}{\sigma} \, \varphi^{(j)}\p{\frac{\mu}{\sigma}},
\end{align*}
where $\varphi^{(j)}$ denotes the $j$-th derivative of the standard Gaussian density with respect to its argument. It is well known that the convolution of two Gaussian random variables is also a Gaussian random variable whose variance is the sum of the original variances, and that convolution commutes with differentiation. In terms of the change of variables $s = x / \sigma$, we see that
\begin{align*}
 &\p{\varphi * \p{\frac{1}{\sigma} \, \varphi\p{\frac{\cdot}{\sigma}} \cdot H_j\p{\frac{\cdot}{\sigma}}}}\p{x} \\
&\ \ \ \ = \frac{1}{\sqrt{1 + \sigma^2}} \frac{1}{\sqrt{j!}} \p{\frac{\partial}{\partial s}}^j \varphi\p{\frac{\sigma}{\sqrt{1 + \sigma^2}} \, s} \\
&\ \ \ \  = \frac{1}{\sqrt{1 + \sigma^2}} \frac{1}{\sqrt{j!}}  \p{\frac{\sigma}{\sqrt{1 + \sigma^2}}}^j \varphi^{(j)}\p{\frac{x}{\sqrt{1 + \sigma^2}}}.
\end{align*}
Plugging this expression into our previous formula, we find that
\begin{align*}
\int &\nu_j\p{\mu_1} P_{g} \p{\mu_1, \, \mu_2} \nu_k\p{\mu_2} d\mu_1 d\mu_2  \\
&= \frac{1}{\sqrt{1 + \sigma^2}} \p{\frac{\sigma}{\sqrt{1 + \sigma^2}}}^{j + k} \int \frac{\varphi^{(j)}\p{\frac{x}{\sqrt{1 + \sigma^2}}} \, \varphi^{(k)}\p{\frac{x}{\sqrt{1 + \sigma^2}}}}{\sqrt{j!} \, \sqrt{k!} \ \varphi\p{\frac{x}{\sqrt{1 + \sigma^2}}}} dx \\
&=   \p{\frac{\sigma}{\sqrt{1 + \sigma^2}}}^{j + k}  \int H_j\p{\frac{x}{\sqrt{1 + \sigma^2}}} H_k\p{\frac{x}{\sqrt{1 + \sigma^2}}} \\
&\ \ \ \ \ \ \ \ \ \ \ \  \frac{1}{\sqrt{1 + \sigma^2}} \, \varphi\p{\frac{x}{\sqrt{1 + \sigma^2}}} dx \\
&=   \p{\frac{\sigma}{\sqrt{1 + \sigma^2}}}^{j + k}  \int H_j\p{x} H_k\p{x}  \varphi\p{x} dx \\
&=  \p{\frac{\sigma}{\sqrt{1 + \sigma^2}}}^{j + k} \delta\p{\cb{j = k}}, 
\end{align*}
because the Hermite polynomials as defined in \eqref{eq:hermite} are orthonormal with respect to the standard Gaussian distribution. By Theorem \ref{theo:most_favorable}, we thus conclude that the most favorable family is given by the first $p$ Hermite polynomials. Moreover, again by Theorem \ref{theo:most_favorable}, the relative efficiency coefficient of this family corresponds to the $p+1$-st eigenvalue of $P_g$, i.e., $ \p{{ \sigma^2}/\p{1 + \sigma^2}}^j$. Because $\Gamma_j$ was a most favorable family for density deconvolution, any other family of the form \eqref{eq:expfam} will have worse relative efficiency.

\subsection*{Proof of Lemma \ref{lemm:pinsker}}

Continuing the argument our argument from Section \ref{sec:gauss_sequence}, it remains to derive a lower bound for the minimax risk $R^L$ among linear estimators in the Gaussian sequence model. We begin by noting that
$$\kappa^{J_C} \leq \mu_C < \kappa^{J_C + 1}, \ \ J_C = \left\lfloor \log\p{\mu_C} \, \big/ \, \log\p{\kappa} \right\rfloor. $$
Thus, we can expand out \eqref{eq:implicit} as
\begin{align*}
C^2 &= \sum_{j = 1}^{J_C} r^{2j}_\sigma \kappa^j \p{\mu_C - \kappa^j} \\
&= \mu_C \, \frac{\p{r_\sigma^2 \kappa}^{{J_C} + 1} - r_\sigma^2\kappa}{r_\sigma^2\kappa - 1} - \frac{\p{r_\sigma^2 \kappa^2}^{{J_C} + 1} - r_\sigma^2\kappa^2}{r_\sigma^2\kappa^2 - 1} \\
&= \p{r_\sigma^2 \kappa^2}^{{J_C} + 1} \, \p{\frac{\mu_C \, \kappa^{-\p{{J_C} + 1}}}{r_\sigma^2\kappa - 1} - \frac{1}{r_\sigma^2\kappa^2 - 1}} - \p{\frac{\mu_C \, r_\sigma^2\kappa}{r_\sigma^2\kappa - 1} - \frac{r_\sigma^2\kappa^2}{r_\sigma^2\kappa^2 - 1}} \\
&\leq B^2_{\sigma, \, \kappa} \p{r_\sigma^2 \kappa^2}^{{J_C} + 1} - 1, \ \ B^2_{\sigma, \, \kappa} = \frac{r_\sigma^2 \kappa \p{\kappa - 1}}{\p{r_\sigma^2 \kappa^2 - 1}\p{r^2 \kappa - 1}},
\end{align*}
where the inequality on the third line hold whenever $C$ (and thus also $\mu_C$ and $J_C$) are large enough. Thus, we conclude that
$$ J_C + 1 \geq \frac{\log\p{\sqrt{C^2 + 1} / B_{\sigma, \, \kappa}}}{\log\p{r_\sigma \kappa}}. $$
We can also bound the risk $R^L$ in \eqref{eq:rl} by
\begin{align*}
R^L
&= \sum_{j = 1}^{J_C} r_\sigma^{2j} \p{1 - \frac{\kappa^j}{\mu_C}} \\
&\geq \sum_{j = 1}^{J_C} r_\sigma^{2j} \p{1 - \frac{\kappa^j}{\kappa^{J_C + 1}}} \\
&= \frac{r_\sigma^{2\p{J_C + 1}} - r_\sigma^2}{r_\sigma^2 - 1} - \frac{\p{r^2_\sigma \kappa}^{J_C + 1} - r_\sigma^2 \kappa}{ \kappa^{J_C + 1} \p{r_\sigma^2 \kappa - 1}} \\
&\geq r_\sigma^{2\p{J_C + 1}} \, \frac{r_\sigma^2 \p{\kappa - 1}}{\p{r_\sigma^2 - 1}\p{r_\sigma^2 \kappa - 1}} - \frac{r_\sigma^2}{r_\sigma^2 - 1}.
\end{align*}
 Plugging in our previous bound for $J_C$, we find that
 \begin{align*}
 R_L &\geq \exp\sqb{\log\p{\frac{\sqrt{C^2 + 1}}{B_{\sigma, \, \kappa}}} \, \frac{2 \log\p{r_\sigma}}{\log\p{r_\sigma} + \log\p{\kappa}}} \, \frac{r_\sigma^2 \p{\kappa - 1}}{\p{r_\sigma^2 - 1}\p{r_\sigma^2 \kappa - 1}}  - \frac{r_\sigma^2}{r_\sigma^2 - 1}\\
&\geq \exp\sqb{\log\p{\frac{C}{B_{\sigma, \, \kappa}}} \, \frac{2 \log\p{r_\sigma}}{\log\p{r_\sigma} + \log\p{\kappa}}} \, \frac{r_\sigma^2 \p{\kappa - 1}}{\p{r_\sigma^2 - 1}\p{r_\sigma^2 \kappa - 1}} ,
 \end{align*}
where the last inequality again holds for large enough $C$.
Once paired with the 5/4 bound of \citet{donoho1990minimax} for Pinsker's constant, this bound yields the desired result.

\subsection*{Proof of Lemma \ref{lemm:poly_risk}}

By the same argument as in the proof of \ref{lemm:gauss_sequence}, we can verify that
\begin{align*}
\lim_{M \rightarrow \infty} \limn \EE{L_n\p{\hgnp}}
&= \limn \sum_{j = 1}^\infty \EE{Z_j - \gamma_j}^2 \\
&= \sum_{j = 1}^p r_\sigma^{2j} + \sum_{j = p + 1}^\infty \gamma_j^2.
\end{align*}
The above expression is largest if the signal concentrates in the $p+1$-st Hermite coefficient, i.e., $\tau\p{\mu} = \gamma_{p+1} H_j\p{\mu/\sigma}$, yielding
\begin{align*}
\limsup_{n \rightarrow \infty} \cb{ \sup_{\tau \in A^\sigma_{\kappa, \, C}}  {L_n\p{\hgnp}} }
&=  \sum_{j = 1}^p r_\sigma^{2j} + C^2 \kappa^{-2\p{p + 1}} \\
&\leq \frac{1}{r_\sigma^2 - 1} r_\sigma^{2\p{p + 1}} + C^2 \kappa^{-2\p{p + 1}}.
\end{align*}
Plugging in the choice for $p$ specified in \eqref{eq:pchoice} and assuming that $C$ is large enough that $\lceil\log(C/\sigma)/\log(r_\sigma\kappa)\rceil - 1 \geq 2$, we get that
\begin{align*} \\
\limsup_{n \rightarrow \infty} \cb{ \sup_{\tau \in A^\sigma_{\kappa, \, C}}  {L_n\p{\hgnp}} } 
&\leq \sigma^2 r^{2\p{\frac{\log\p{C / \sigma}}{\log{r_\sigma \kappa}} + 1}} + C^2 \kappa^{-2\p{\frac{\log\p{C / \sigma}}{\log{r_\sigma \kappa}}}} \\
&= \p{1 + r^2} \sigma^\frac{2\log\p{k}}{\log\p{r_\sigma k}} C^\frac{2\log\p{r_\sigma}}{\log\p{r_\sigma k}},
\end{align*}
which is what we set out to show.

\subsection*{Proof of Lemma \ref{lemm:gauss_sequence}}

Our proof proceeds in several parts. We begin by establishing a version of our result for a simpler finite-dimensional problem; the general statement then shows that the finite-dimensional case can uniformly approximate our problem of interest.

\paragraph{A Finite-Dimensional Model} For some $J \in \NN$, suppose that $\tau$ is known to lie in an ellipse $ A^\sigma_{\kappa, \, C}\p{J}$ defined by
\begin{equation*}
\label{eq:ellipse_finite}
\tau\p{\mu} = \sum_{j = 1}^J \gamma_j \, H_j\p{\frac{\mu}{\sigma}}, \ \ \gamma \in \ell_2^\kappa\p{C, \, J} := \cb{\gamma' : \sum_{j = 1}^J \kappa^{2j} \gamma_j^2 \leq C^2},
\end{equation*}
and that we only consider estimators over the set $\Lambda_\kappa^\sigma\p{J}$ defined by
\begin{equation*}
\hgn\p{\mu} = g^M_\sigma\p{\mu} \exp\sqb{\frac{1}{\sqrt{n}} \sum_{j = 1}^J \hgamma_j \, H_j\p{\frac{\mu}{\sigma}} - \psi_n\p{\hgamma}}, \ \  \sum_{j = 1}^J \kappa^{2j} \hgamma_j^2 < L^2.
\end{equation*}
Our first task is to show that the the minimax risk over this finite-dimensional parametric class can, for large $M$, be well-approximated by minimax risk of the analogous finite Gaussian problem. We recall that, as usual, $\smash{g^M_\sigma\p{\mu}}$ denotes the Gaussian density $g_\sigma\p{\mu}$ that has been ``wrapped around'' the interval $\Omega_M = [-M, \, M]$.

\paragraph{Convergence of the Likelihood}

Consider any parameter $\gamma'$ whose induced tilting function satisfies $\smash{\tau_{\gamma'} \in A^\sigma_{\kappa, \, C}\p{J}}$, and denote the resulting marginal density function by $\smash{f^{(n)}_{\gamma'} \propto \varphi * g^M_\sigma \, e^{\tau_{\gamma'} / \sqrt{n}}}$. Because the basis functions $H_j$ are all bounded (recall that we assume a compact domain $[-M, \, M]$), the log-likelihood $f^{(n)}_{\gamma'}\p{x}$ is uniformly Lipschitz in $\gamma'$ for all $x$. Thus, we can use standard empirical process theory results to verify that the log-likelihood at $\gamma'$, namely $\smash{\log(\prod_{i = 1}^n f^{(n)}_{\gamma'}(X_i))}$, is entirely determined by the score at the optimum $\gamma$ \citep[e.g.,][Lemma 19.31]{van2000asymptotic}:
\begin{equation*}
\sup_{\tau_{\gamma'} \in A^\sigma_{\kappa, \, C}\p{J}} \abs{\sqrt{n} \log\p{\frac{\prod_{i = 1}^n f^{(n)}_{\gamma'}\p{X_i}}{\prod_{i = 1}^n f^{(n)}_{\gamma}\p{X_i}}} -  \p{\gamma' - \gamma} \cdot \sum_{i = 1}^n \nabla \log \p{f^{(n)}_{\gamma}\p{X_i}}} \rightarrow_p 0.
\end{equation*}
As $\gamma$ is the true optimal parameter, we know that $\smash{\EE[\gamma]{\nabla \log \p{f^{(n)}_{\gamma}\p{X_i}}} = 0}$; meanwhile, $\smash{n \Var[\gamma]{\nabla \log \p{f^{(n)}_{\gamma}\p{X_i}}}}$ converges to the Fisher information at $\gamma = 0$.
Thus, by the central limit theorem, we can verify that
\begin{equation*}
\sum_{i = 1}^n \nabla \log \p{f^{(n)}_{\gamma}\p{X_i}} \Rightarrow \nn\p{0, \, V_{M, \,J}}, \ \ V_{M, \,J} := \Var[0]{\nabla \log \p{f_{0}\p{X_i}}};
\end{equation*}
Thus, we conclude that the log-likelihood in favor of $\gamma'$ relative to $\gamma$ is asymptotically equivalent in distribution to the log-likelihood arising from the Gaussian experiment where we observe
$\smash{Z_M \sim  \nn\p{\gamma, \, V_{M, \,J}}}$
and want to recover $\smash{\gamma \in \ell_2^\kappa\p{C, \, J}}$.

\paragraph{Convergence of the Loss}

Similarly, we can verify that the density estimation loss at $\gamma'$ given the true parameter value $\gamma$ satisfies
\begin{align*}
&L_n\p{g^{(n)}_{\tau_{\gamma'}}} 
 = n \int_{\Omega_M} g^M_\sigma\p{\mu} e^{\frac{1}{\sqrt{n}} \gamma \cdot H_{(1:J)} \p{\frac{\mu}{\sigma}} - \psi_M\p{\frac{\gamma}{\sqrt{n}}}}  \\
  &\ \ \ \ \ \ \ \ \p{\frac{1}{\sqrt{n}} \p{\gamma - \gamma'} \cdot  H_{(1:J)} \p{\frac{\mu}{\sigma}} -\psi_M\p{\frac{\gamma}{\sqrt{n}}} + \psi_M\p{\frac{\gamma'}{\sqrt{n}}} } \ d\mu, \\
 \end{align*}
 and so
 \begin{align*}
 \limn L_n \p{g^{(n)}_{\tau_{\gamma'}}} 
 &=  \int_{\Omega_M} \p{\p{\gamma - \gamma'} \cdot \p{H_{(1:J)}\p{\frac{\mu}{\sigma}} - \psi'_M\p{0}}}^2 \,  g^M_\sigma\p{\mu}  \ d\mu \\
 &= \p{\gamma' - \gamma}^\top W_{M, \,J} \p{\gamma' - \gamma},\\
  & \hspace{-10mm}  \with W_{M, \,J} := \int_{\Omega_M} \p{H_{(1:J)}\p{\frac{\mu}{\sigma}} - \psi'_M\p{0}}^{\otimes 2} \,  g^M_\sigma\p{\mu}  \ d\mu,
\end{align*}
where $\smash{H_{(1:J)}(\cdot)}$ is a vector obtained by stacking the first $J$ Hermite functions. Moreover, this convergence is uniform in $\gamma'$. Thus, for large $n$, our re-scaled deviance loss is asymptotically equal to
$$ L_{M, \,J}\p{\gamma'} := \p{\gamma' - \gamma}^\top W_{M, \,J} \p{\gamma' - \gamma}. $$

\paragraph{Convergence of the Minimax Risk}

Given the convergence results for the log-likelihood and for the loss established above, we should expect the statistical problem of density estimation to be asymptotically equivalent to finding the mean of a Gaussian vector with variance $\smash{V_{M, \,J}}$ under loss $\smash{L_{M, \,J}}$. To establish this formally, we can use least-favorable priors. Let $\smash{\hgamma^*_{M, \, J}}$ be a minimax estimator for $\gamma$ in the Gaussian model with variance $\smash{V_{M, \,J}}$ under loss $\smash{L_{M, \,J}}$, subject to  $\gamma \in \ell_2^\kappa\p{C, \, J}$. By checking the conditions of \citet{wald1945statistical}, we can verify that $\hgamma^*_{M, \, J}$ is unique and that it is also Bayes for a least-favorable prior $\pi^*_{M, \, J}$; moreover, the risk of $\hgamma^*_{M, \, J}$ is constant over the support of $\pi^*_{M, \, J}$. Finally, because $L_{M, \,J}$ is quadratic, $\hgamma^*_{M, \, J}$ is the posterior mean for this least-favorable prior.

Now, let $\hgamma^{(n, \, \pi^*)}_{M, \, J}$ be the posterior mean for $\gamma$ in the density estimation problem with $n$ samples, where $\gamma$ has a prior $\pi^*_{M, \, J}$. By our previous results, we find that
$$ \law\p{L_n\p{g^{(n)}_{\tau_{\hgamma^(n, \, \pi^*)_{M, \, J}}}}}
\Rightarrow \law\p{L_{M, \, J} \p{\hgamma^*_{M, \, J}}}, $$
and that the risk of this estimator is asymptotically constant over the support of $\pi^*_{M, \, J}$; thus, this estimator is asymptotically minimax over the support of $\pi^*_{M, \, J}$. Finally, the risk of this estimator is never asymptotically worse outside the support of $\pi^*_{M, \, J}$, and so $\smash{\hgamma^{(n, \, \pi^*)}_{M, \, J}}$ is in fact asymptotically minimax for our finite-dimensional density estimation problem.

\paragraph{Taking Limits}

To move from the compact interval $\Omega_M$ to the real line we observe that, for any fixed $J$, we know that the Hermite functions are orthogonal over the whole real line:
\begin{align*}
&\int_{\RR} H_{(1:J)}\p{\frac{\mu}{\sigma}} \,  g_\sigma\p{\mu}  \ d\mu = 0, \\
&\int_{\RR} H_{(1:J)}\p{\frac{\mu}{\sigma}}^{\otimes 2} \,  g_\sigma\p{\mu}  \ d\mu = I_{J \times J},
\end{align*}
and so $\lim_{M \rightarrow \infty} W_{M, \, J} = I_{J \times J}$. Moreover, by the same argument as used in the proof of Theorem \ref{theo:gauss_hard}, we find that $V_{M, \, J}$ converges to the leading $J \times J$ sub-matrix of the covariance matrix defined in \eqref{eq:gauss_sequence}. Meanwhile, for a fixed $M$, we can use the uniform approximability result from Corollary \ref{coro:approx} to verify that the minimax risk of density estimation converges as $J \rightarrow \infty$; a similar argument holds in the Gaussian case. Because the minimax risks for both problems match for every finite $J$, they must thus also match as $J \rightarrow \infty$. Now, because we have established convergence both when $J$ goes to infinity given a fixed $M$ and when $M$ goes to infinity given a fixed $J$, we conclude that the joint limit $M, \, J \rightarrow \infty$ is well-defined and does not depend on the order in which we take the limits; this implies the desired result.

\subsection*{Proof of Lemma \ref{lemm:general_coeff}}

Asymptotic normality follows directly from the central limit theorem; we only need to verify the moments. Now, using similar arguments as in the proof of Corollary \ref{coro:approx}, we can verify that the limiting covariance of the $Z_j^{(n)}$ does not depend on $\tau$, and that
\begin{align*}
\limn & \Cov{Z_j^{(n)}, \, Z_{j'}^{(n)}} =  \rho_j^{-2} \,  \int_\Omega U_j(x) U_{j'}(x) f_0(x) \ dx \\
&=  \rho_j^{-2} \int_{\Omega^3}  \frac{K\p{x, \, \mu_1} T_j\p{\mu_1} g_0\p{\mu_1} \,  K\p{x, \, \mu_2} T_{j'}\p{\mu_2} g_0\p{\mu_2} }{f_0\p{x}} \ d\mu_1 \ d\mu_2 \ dx \\
&= \rho_j^{-2}  \int_{\Omega^2} \zeta_{j+1}\p{\mu_1} \, P_{g_0}\p{\mu_1, \, \mu_2} \, \zeta_{j' + 1}\p{\mu_2}  \ d\mu_1 \ d\mu_2 \\
&=  \delta_{\cb{j = j'}} \, \rho_j^{-1},
\end{align*}
where $P_{g_0}$ is the linear operator defined in \eqref{eq:lin} and the $\zeta_j$ are its eigenvectors as defined in Theorem \ref{theo:most_favorable}; recall that $\rho_j = \lambda_{j+1}\p{P_{g_0}}$. Meanwhile, the $T_j$ are centered such that $\smash{\EE{T_j(\mu)} = 0}$, and so
\begin{align*}
\frac{1}{\sqrt{n}} &\EE[0]{Z_j^{(n)}}
=  \int_\Omega U_j(x) f_0(x) \ dx  \\
&=  \int_{\Omega} K\p{x, \, \mu} \int_\Omega T_j\p{\mu} g_0\p{\mu} \ d\mu \ dx 
= 0.
\end{align*}
Thus, we can verify that
\begin{align*}
& \limn \EE{Z_j^{(n)}} =  \limn \sqrt{n}  \,  \rho_j^{-1} \, \int_\Omega U_j(x) \, f_{\tau / \sqrt{n}}(x) \ dx \\
&\ \ =\rho_j^{-1} \int_\Omega U_j(x) \sqb{\frac{\partial}{\partial \varepsilon} f_{\varepsilon \tau} (x)}_{\varepsilon = 0} \ dx \\
&\ \ = \rho_j^{-1} \int_{\Omega^3}  \frac{K\p{x, \, \mu_1} T_j\p{\mu_1} g_0\p{\mu_1} \,  K\p{x, \, \mu_2} \sqb{\frac{\partial}{\partial \varepsilon} g_{\varepsilon \tau} \p{\mu_2}}_{\varepsilon = 0}}{f_0\p{x}} \ d\mu_1 \ d\mu_2 \ dx \\
&\ \ = \rho_j^{-1} \int_{\Omega^2} T_j \p{\mu_1} \, \sqrt{g_0\p{\mu_1}} \, P_{g_0} \, \sqrt{g_0\p{\mu_2}} \,  \tau\p{\mu_2}  \ d\mu_1 \ d\mu_2 \\
&\ \ = \int_\Omega T_j\p{\mu} \, \tau\p{\mu} \, g_0\p{\mu} \ d\mu,
\end{align*}
where the last equality follow from the spectral theorem because $\smash{T_j \sqrt{g_0}}$ is an eigenfunction of $\smash{P_{g_0}}$ with eigenvalue $\rho_j$.

\end{appendix}

\end{document}